\documentclass{amsart}
\usepackage{amssymb}
\usepackage{amsfonts}

\newtheorem{theorem}{Theorem}
\theoremstyle{plain}

\newtheorem{corollary}[theorem]{Corollary}

\newtheorem{definition}[theorem]{Definition}

\newtheorem{lemma}[theorem]{Lemma}

\newtheorem{proposition}[theorem]{Proposition}
\newtheorem{remark}[theorem]{Remark}

\numberwithin{theorem}{section}
\input{tcilatex}

\begin{document}
\title[Reduction $\func{mod}p$ of Cuspidal Representations of $GL_{2}(%
\mathbb{F}_{p^{n}})$]{Reduction $\func{mod}p$ of Cuspidal Representations of 
$GL_{2}(\mathbb{F}_{p^{n}})$ and Symmetric Powers}
\author{Davide A. Reduzzi}
\address{University of California at Los Angeles\\
520 Portola Plaza, Math Sciences Building 6363\\
Los Angeles, CA 90095, USA}

\begin{abstract}
We show the existence of integral models for cuspidal representations of $%
GL_{2}(\mathbb{F}_{q})$, whose reduction modulo $p$ can be identified with
the cokernel of a differential operator on $\mathbb{F}_{q}[X,Y]$ defined by
J-P. Serre. These integral models come from the crystalline cohomology of
the projective curve $XY^{q}-X^{q}Y-Z^{q+1}=0$.

As an application, we can extend a construction of C. Khare and B. Edixhoven
(\cite{K-E})\ giving a cohomological analogue of the Hasse invariant
operator acting on spaces of $\func{mod}p$ modular forms for $GL_{2}$.
\end{abstract}

\email{devredu83@math.ucla.edu}
\maketitle

\section{Introduction}

Let $q=p^{n}$ be a power of a prime number $p$; the irreducible complex
representations of $G=GL_{2}(\mathbb{F}_{q})$ (of dimension greater than
one) that are not twists of the Steinberg representation are of two types:\
the principal series representations, having dimension $q+1$ and obtained by
inducing to $G$ characters of the Borel subgroup of $G$, and the cuspidal
representations, having dimension $q-1$ and characterized by the property
that they do not occur as a factor of a principal series.

The dimensions of the cuspidal representations and of the principal series
representations of $G$ appear in the study of the modular representations of
the group: for $k\geq 0$ let $V_{k}=\limfunc{Sym}^{k}\mathbb{F}_{q}^{2}$ be
the $k^{th}$-symmetric power of the standard left representation of $G$ over
the field $\mathbb{F}_{q}$; in \cite{Serre}, J-P. Serre noted the following
identity in the Grothendieck group of $G$ (cf. Proposition \ref{identit}):%
\begin{equation*}
V_{k}-e\cdot V_{k-(q+1)}=V_{k-\left( q-1\right) }-e\cdot V_{k-2q}\text{ \ \
\ \ \ (}k\in 
\mathbb{Z}
\text{),}
\end{equation*}

\noindent where the definition of the $V_{k}$'s (as elements of the
Grothendieck group of $G$) has been \textit{extended} in a suitable way for $%
k<0$ (cf. Definition \ref{newww}), and $e$ denotes the character determinant.

A look at the left hand side of the above relation leads naturally to
consider the $\mathbb{F}_{q}[G]$-equivariant\ map:%
\begin{equation*}
\theta _{q}:e\otimes V_{k-(q+1)}\rightarrow V_{k}
\end{equation*}

\noindent defined, for any $k>q$, as multiplication by the Dickson invariant 
$X^{q}Y-XY^{q}\in \mathbb{F}_{q}[X,Y]$ (cf. \cite{Glo}) and realizing $%
e\otimes V_{k-(q+1)}$ as a subobject of $V_{k}$. The cokernel of this map
has dimension $q+1$ and turns out to be the reduction $\func{mod}p$ of a
principal series representation of $G$: $\limfunc{coker}\theta _{q}\simeq 
\limfunc{Ind}_{B}^{G}\left( \eta ^{k}\right) $, where $B\leq G$ is the
subgroup of upper triangular matrices, and $\eta $ is the character of $B$
given by $\left( 
\begin{array}{cc}
a & b \\ 
0 & c%
\end{array}%
\right) \mapsto a$. In particular, one deduces the following periodicity
result (cf. Proposition \ref{fame}):

\bigskip

\noindent \textbf{Proposition.} \textit{For any integer }$k>q$\textit{\ and
any }$\lambda \geq 0$\textit{\ there are isomorphisms of }$G$\textit{%
-modules:}%
\begin{equation*}
\frac{V_{k}}{e\otimes V_{k-(q+1)}}\simeq \limfunc{Ind}\nolimits_{B}^{G}%
\left( \eta ^{k}\right) \simeq \frac{V_{k+\lambda (q-1)}}{e\otimes
V_{k+\lambda (q-1)-(q+1)}}.
\end{equation*}

\bigskip

It seems natural to expect that the other "period" $q-1$ of the group $G$ is
associated to reduction $\func{mod}p$ of cuspidal representations of $G$.
\noindent At this purpose, in response to a message of C. Khare, Serre
defined a $G$-map $D:V_{k}\rightarrow V_{k+(q-1)}$ that raises the power of
the symmetric module $V_{k}$ by $q-1$. This function is defined as a
derivation map of the algebra $\mathbb{F}_{q}[X,Y]$, as follows:%
\begin{equation*}
D:\mathbb{F}_{q}[X,Y]\rightarrow \mathbb{F}_{q}[X,Y]:\text{ \ }f\longmapsto
X^{q}\frac{\partial f}{\partial X}+Y^{q}\frac{\partial f}{\partial Y}.
\end{equation*}

The derivation map $D$ in not always injective (cf. Proposition \ref%
{injective}) and seems to capture essential properties related to the
existence or non-existence of embeddings of $G$-modules of the form $%
V_{k}\rightarrow V_{k+(q-1)}$ (cf. Proposition \ref{essential}).

In this paper we show how the cokernel of the $D$-map can be identified with
the reduction of a cuspidal representation of $G$ over the field $\overline{%
\mathbb{Q}
}_{p}$; more precisely, in sections $3$ and $4$ we consider, for $2\leq
k\leq p-1$, the exact sequence of $G$-modules (Proposition \ref{sonno2}):%
\begin{equation*}
0\rightarrow e\otimes V_{k-2}\overset{\overline{\theta }_{q}}{\rightarrow }%
\dfrac{V_{k+(q-1)}}{D(V_{k})}\rightarrow \limfunc{coker}\overline{\theta }%
_{q}\rightarrow 0
\end{equation*}

\noindent and we identify it with an exact sequence constructed by B.
Haastert and J.C. Jantzen in \cite{HaaJan}, and coming from the crystalline
cohomology of the Deligne-Lusztig variety of the group $SL_{2/\mathbb{F}%
_{q}} $, i.e. the smooth projective curve $C_{/\mathbb{F}%
_{q}}:XY^{q}-X^{q}Y-Z^{q+1}=0$. We then deduce the main result of the paper,
analogous to the periodicity result for reduction of principal series
representations of $G$ (cf. Theorem \ref{fame2}):

\bigskip

\noindent \textbf{Theorem. }\textit{Let }$q>2$\textit{, }$2\leq k\leq p-1$%
\textit{\ with }$k\neq \frac{q+1}{2}$\textit{\ and let us denote by }$\Theta
\left( \chi ^{k}\right) $\textit{\ the cuspidal }$\overline{%
\mathbb{Q}
}_{p}$\textit{-representation of }$G$\textit{\ associated to the }$k^{th}$%
\textit{-power of the Teichm\"{u}ller character }$\chi $\textit{. Then there
exists a natural integral model }$\widetilde{\Theta }\left( \chi ^{k}\right) 
$\textit{\ of }$\Theta \left( \chi ^{k}\right) $\textit{\ coming from the }$%
-k$\textit{-eigenspace of the first crystalline cohomology group of the
projective curve }$C_{/\mathbb{F}_{q}}$\textit{\ such that there is an
isomorphism of }$\mathbb{F}_{q}\left[ G\right] $\textit{-modules:}%
\begin{equation*}
\dfrac{V_{k+(q-1)}}{D\left( V_{k}\right) }\simeq \widetilde{\Theta }\left(
\chi ^{k}\right) \text{ }(\func{mod}p).
\end{equation*}

\bigskip

In the last section of the paper, we switch to modular forms $\func{mod}p$
for $GL_{2}$: in this context, thanks to the Eichler-Shimura isomorphism,
the modules $V_{k}$'s appear in the study of the systems of eigenforms for
the action of the Hecke algebra on the space of modular forms $M_{k}\left(
\Gamma _{1}(N),\overline{\mathbb{F}}_{p}\right) $. As a consequence, the
maps $\theta _{p}$ and $D$ considered above play a role too: in \cite{A-S},
A. Ash and G. Stevens identify a group theoretical analogue of the theta
operator in the map induced in cohomology by the Dickson homomorphism $%
\theta _{p}$; in \cite{K}, B. Edixhoven and C. Khare build a cohomological
analogue of the Hasse invariant operator as a map $\alpha :H^{1}(\Gamma
_{1}(N),V_{0})\hookrightarrow H^{1}(\Gamma _{1}(N),V_{p-1})$. The existence
of the map $D$ might lead to further results in this direction.

\bigskip \newpage 

\textbf{Acknowledgements}

I am deeply grateful to Chandrashekhar Khare for conceiving the problems
studied here: the viewpoint on cohomological weight shiftings appearing on
this paper is entirely his. I further thank him for the generous help he
gave me in overcoming crucial points of the proofs appearing in the paper,
and for sharing his correspondence with Gordan Savin and Jean-Pierre Serre
on the topics considered here. The question answered by Theorem \ref{fame2}
is inspired by J-P.Serre's letter to J. Tate (\cite{Se2}), and the
derivation map $D$ that J-P. Serre defined in a correspondence with C. Khare.

I thank Najmuddin Fakhruddin, for signaling that the reference \cite{HaaJan}
would help in proving Theorem \ref{fame2}, and Gordan Savin, for an
interesting conversation we had. I would like to further express my
gratitude to Florian Herzig, for reading this paper and giving me
suggestions for its improvement.

\section{Reduction of principal series representations}

We introduce some notation that we will be using throughout the paper: let $%
q=p^{n}>1$ be an integral power of a prime number $p$ and let $\mathbb{F}=%
\mathbb{F}_{q}$ be a field of cardinality $q$, with algebraic closure $%
\overline{\mathbb{F}}$; $G=GL_{2}\left( \mathbb{F}\right) $ will denote the
group of $2\times 2$ non-singular matrices with coefficients in $\mathbb{F}$%
, and $V=$ $\mathbb{F}^{2}$ the standard left $G$-representation.

\subsection{Definition of the $G$-modules $V_{k}$}

For any non-negative integer $k$ let $V_{k}=\limfunc{Sym}^{k}V$ be the
symmetric $k^{th}$-power of the $G$-representation $V$: $V_{k}$ will be
identified with the $\mathbb{F}$-vector space $\mathbb{F}\left[ X,Y\right]
_{k}$ of homogeneous polynomials over $\mathbb{F}$ in two variables and of
degree $k$, endowed with the left action of $G$ induced by:

\begin{center}
\begin{equation}
\left( 
\begin{array}{cc}
a & b \\ 
c & d%
\end{array}%
\right) \cdot X=aX+cY,\ \left( 
\begin{array}{cc}
a & b \\ 
c & d%
\end{array}%
\right) \cdot Y=bX+dY.  \label{action}
\end{equation}
\end{center}

\noindent It is possible to extend the definition of the $V_{k}$'s - as
elements of the Grothendieck group of $G$ - to the negative integers: this
is done by considering (a variant of) the Euler-Poincar\'{e} characteristic
of the twists $\mathcal{O}\left( k\right) $ of the structure sheaf on $%
\mathbb{P}_{/\mathbb{F}}^{1}$ (cf. Serre, \cite{Serre}). In order to justify
the presence of this sheaf, we will start by considering the case of the
algebraic group $GL_{2}$.

Let $\mathbf{G}=GL_{2}$ as an algebraic group over $\mathbb{F}$, and let $%
\mathbf{T\subset G\ }$be the maximal split torus of diagonal matrices. We
identify the character group $X(\mathbf{T})$ of $\mathbf{T}$ with $%
\mathbb{Z}
^{2}$ in the usual way, so that the roots associated to the pair $\left( 
\mathbf{G},\mathbf{T}\right) $ are $(1,-1)$ and $(-1,1)$;\ we fix a choice
of positive root $\alpha =(1,-1)$. The corresponding Borel subgroup $\mathbf{%
B}$ is the group of upper triangular matrices in $\mathbf{G}$; we denote by $%
\mathbf{B}^{-}$ the opposite Borel subgroup.

For a fixed $\lambda \in $ $X(\mathbf{T})$, let $M_{\lambda }$ be the one
dimensional left $\mathbf{B}^{-}$-module on which $\mathbf{B}^{-}$ acts
(through its quotient $\mathbf{T}$) via the character $\lambda $ . Denote by 
$\limfunc{ind}\nolimits_{\mathbf{B}^{-}}^{\mathbf{G}}M_{\lambda }$ the left $%
\mathbf{G}$-module given by algebraic induction from $\mathbf{B}^{-}$ to $%
\mathbf{G}$ of $M_{\lambda }:$%
\begin{equation*}
\limfunc{ind}\nolimits_{\mathbf{B}^{-}}^{\mathbf{G}}M_{\lambda }:=\left\{
f\in \limfunc{Mor}\left( \mathbf{G},\mathbb{G}_{a}\right) :f\left( bg\right)
=\lambda \left( b\right) f\left( g\right) \text{, for all }g\in \mathbf{G}%
,b\in \mathbf{B}^{-}\right\} ,
\end{equation*}

\noindent where $\limfunc{Mor}\left( \mathbf{B}^{-},\mathbb{G}_{a}\right) $
denotes the set of morphisms of functors from $\mathbf{B}^{-}$ to $\mathbb{G}%
_{a}$, and the left action of $\mathbf{G}$ on the above space is given by
letting $\left( gf\right) \left( x\right) =f\left( xg\right) $ for all $%
x,g\in \mathbf{G}$ and $f\in \limfunc{ind}\nolimits_{\mathbf{B}^{-}}^{%
\mathbf{G}}M_{\lambda }$.

Define the following generalization of the dual \textit{Weyl module} for $%
\lambda $ (cf. \cite{Jan}, II.5):%
\begin{equation*}
W\left( \lambda \right) :=\sum\nolimits_{i\geq 0}\left( -1\right) ^{i}\cdot
R^{i}\limfunc{ind}\nolimits_{\mathbf{B}^{-}}^{\mathbf{G}}\left( M_{\lambda
}\right) \text{,}
\end{equation*}

\noindent where $R^{i}\limfunc{ind}\nolimits_{\mathbf{B}^{-}}^{\mathbf{G}}$
denote the $i^{th}$-right derived functor of the functor $\limfunc{ind}%
\nolimits_{\mathbf{B}^{-}}^{\mathbf{G}}$.

Notice that $W\left( \lambda \right) $ is a well defined element of the
Grothendieck group $K_{0}(\mathbf{G})$ of $\mathbf{G}$, because each $R^{i}%
\limfunc{ind}\nolimits_{\mathbf{B}^{-}}^{\mathbf{G}}\left( M_{\lambda
}\right) $ is a finite dimensional $\mathbf{G}$-module, and $R^{i}\limfunc{%
ind}\nolimits_{\mathbf{B}^{-}}^{\mathbf{G}}\left( M_{\lambda }\right) $ is
zero for $i>1$ (\cite{Jan}, II.4.2).

\noindent By \cite{Jan}, I.5.12, we have for each $i\geq 0$ a canonical
isomorphism of $\mathbf{G}$-modules:%
\begin{equation*}
R^{i}\limfunc{ind}\nolimits_{\mathbf{B}^{-}}^{\mathbf{G}}\left( M_{\lambda
}\right) \simeq H^{i}(\mathbf{B}^{-}\backslash \mathbf{G},\mathcal{%
\tciLaplace }(M_{\lambda })),
\end{equation*}

\noindent where $\mathcal{\tciLaplace }(M_{\lambda })$ is the sheaf on the
projective scheme $\mathbf{B}^{-}\backslash \mathbf{G}\simeq \mathbb{P}^{1}$
associated to $M_{\lambda }$, i.e. for an open subset $U\subseteq \mathbb{P}%
^{1}$ we have:%
\begin{equation*}
\mathcal{\tciLaplace }(M_{\lambda })(U):=\left\{ f\in \limfunc{Mor}\left(
\pi ^{-1}U,\mathbb{G}_{a}\right) :f\left( bx\right) =\lambda (b)f\left(
x\right) \text{, for all }x\in \pi ^{-1}U,b\in \mathbf{B}^{-}\right\} ,
\end{equation*}

\noindent where $\pi :\mathbf{G}\rightarrow \mathbf{B}^{-}\backslash \mathbf{%
G}\simeq \mathbb{P}^{1}$ is the canonical projection (the $\mathbf{G}$%
-module structure on $H^{i}(\mathbf{B}^{-}\backslash \mathbf{G},\mathcal{%
\tciLaplace }(M_{\lambda }))$ comes from the fact that $\mathcal{\tciLaplace 
}(M_{\lambda })$ is a $\mathbf{G}$-linearized sheaf,\ cf. Remark in \cite%
{Jan}, I.5.12).

Let us fix $\lambda _{k}\in X\left( \mathbf{T}\right) $ corresponding to the
pair of integers $\left( k,0\right) $ ($k\in 
\mathbb{Z}
$); we have an isomorphism of sheaves on $\mathbb{P}^{1}$ (cf. \cite{Jan}
II.4.3.):%
\begin{equation*}
\mathcal{\tciLaplace }(M_{\lambda _{k}})\simeq \mathcal{O}\left( k\right) .
\end{equation*}

To see this, let us work for simplicity over $\overline{\mathbb{F}}$ ; the
canonical projection $\pi $ is given by: 
\begin{equation*}
\begin{pmatrix}
x_{1} & x_{2} \\ 
x_{3} & x_{4}%
\end{pmatrix}%
\mapsto (x_{1}:x_{2})\in \mathbb{P}^{1}.
\end{equation*}

Denote by $\omega :$ $\mathbf{B}^{-}\rightarrow \overline{\mathbb{F}}%
^{\times }$ the character associated to the weight $\lambda =(1,0)\in
X\left( \mathbf{T}\right) $, and let $N:=\limfunc{Ker}\omega $. Fix an open
subset $U$ of $\mathbb{P}^{1}$;\ by definition of the sheaf $\mathcal{%
\tciLaplace }(M_{\lambda })$, if $f\in \mathcal{\tciLaplace }(M_{\lambda
})(U)$ then $f\left( b\cdot x\right) =\omega \left( b\right) ^{k}f\left(
x\right) $ for all $x\in \pi ^{-1}U,b\in \mathbf{B}^{-}$, and $f$ is left $N$%
-invariant. We make the identification of sets $N\backslash \pi ^{-1}U\simeq
U$ through the map $Nx\mapsto e_{1}\cdot x$, where $e_{1}=(1,0)$ and $x\in
\pi ^{-1}U$. We deduce the existence of an isomorphism of $\mathcal{%
\tciLaplace }(M_{\lambda })(U)$ with the group:%
\begin{equation*}
\left\{ f\in \limfunc{Mor}\left( U,\overline{\mathbb{F}}\right) :f\left(
e_{1}\cdot b\cdot x\right) =\omega \left( b\right) ^{k}f\left( e_{1}\cdot
x\right) \text{, for all }x\in \pi ^{-1}U,b\in \mathbf{B}^{-}\right\} .
\end{equation*}%
Since $e_{1}\cdot b=$ $\omega \left( b\right) e_{1}$, we deduce:%
\begin{eqnarray*}
&&\mathcal{\tciLaplace }(M_{\lambda })(U) \\
&\simeq &\left\{ f\in \limfunc{Mor}\left( U,\overline{\mathbb{F}}\right)
:f\left( av\right) =a^{k}f\left( v\right) \text{, for all }v\in U,a\in 
\overline{\mathbb{F}}\right\} \\
\;\; &\simeq &\mathcal{O}(k)(U).
\end{eqnarray*}

\bigskip

\noindent We conclude the existence of isomorphisms of $\mathbf{G}$-modules $%
H^{i}(\mathbf{B}^{-}\backslash \mathbf{G},\mathcal{\tciLaplace }(M_{\lambda
_{k}}))\simeq H^{i}(\mathbb{P}^{1},\mathcal{O}\left( k\right) )$ for every $%
i\geq 0$.

If $k\geq 0$, $H^{1}(\mathbb{P}^{1},\mathcal{O}\left( k\right) )=0$ so that $%
W\left( \lambda _{k}\right) =H^{0}(\mathbb{P}^{1},\mathcal{O}\left( k\right)
)=\limfunc{Sym}\nolimits^{k}\overline{\mathbb{F}}^{2}$ in $K_{0}\left( 
\mathbf{G}\right) $; if $k<0$ we have $H^{0}(\mathbb{P}^{1},\mathcal{O}%
\left( k\right) )=0$\ and $W\left( \lambda _{k}\right) =-H^{1}(\mathbb{P}%
^{1},\mathcal{O}\left( k\right) )$; the canonical perfect pairing of $%
\mathbf{G}$-modules: 
\begin{equation*}
H^{0}(\mathbb{P}^{1},\mathcal{O}\left( -k-2\right) )\times H^{1}(\mathbb{P}%
^{1},\mathcal{O}\left( k\right) )\rightarrow H^{1}(\mathbb{P}^{1},\mathcal{O}%
\left( -2\right) )\simeq \det\nolimits^{-1}\otimes \mathbb{G}_{a},
\end{equation*}%
(cf. \cite{Har}, III.5) gives the following equalities in $K_{0}\left( 
\mathbf{G}\right) $: 
\begin{eqnarray*}
W\left( \lambda _{k}\right) &=&-\limfunc{Hom}\nolimits_{\mathbf{G}}(H^{0}(%
\mathbb{P}^{1},\mathcal{O}\left( -k-2\right) ),\det\nolimits^{-1}\otimes 
\mathbb{G}_{a})= \\
&=&-\limfunc{Hom}\nolimits_{\mathbf{G}}(\det \otimes H^{0}(\mathbb{P}^{1},%
\mathcal{O}\left( -k-2\right) ),\mathbb{G}_{a})= \\
&=&-\left( \det \otimes H^{0}(\mathbb{P}^{1},\mathcal{O}\left( -k-2\right)
)\right) ^{\ast }= \\
&=&-e^{1+k}\cdot \limfunc{Sym}\nolimits^{-k-2}\overline{\mathbb{F}}^{2},
\end{eqnarray*}%
where $\cdot $ denotes the product induced in $K_{0}\left( \mathbf{G}\right) 
$ by the tensor product of $\mathbf{G}$-modules, and $e$ is the determinant
character.

\begin{remark}
In more general settings, let $\mathbf{G}$ be a connected split reductive
group over $\mathbb{F}$; the vanishing of the modules $R^{i}\limfunc{ind}%
\nolimits_{\mathbf{B}^{-}}^{\mathbf{G}}\left( M_{\lambda }\right) $ for $i>0$
whenever $\lambda $ is a dominant weight for $(\mathbf{G,T})$ is the content
of Kempf's vanishing Theorem (cf \cite{Jan}, II.4); on the other side $%
\lambda \in X\left( \mathbf{T}\right) $ is dominant if and only if $\limfunc{%
ind}\nolimits_{\mathbf{B}^{-}}^{\mathbf{G}}\left( M_{\lambda }\right) \neq 0$
(cf. \cite{Jan}, II.2.6). In our case, with the choice of positive root $%
\alpha $, $\lambda _{k}$ is dominant if and only if $k\geq 0$.
\end{remark}

\bigskip

The above considerations on the generalized Weyl modules $W\left( \lambda
_{k}\right) $ for $\mathbf{G}$ lead to the following definition for the
finite group $G$:

\begin{definition}
\label{newww} For any integer $k$, define the element $V_{k}$ of the
Grothendieck group $K_{0}\left( G\right) $ of $G$ over $\mathbb{F}$ by:%
\begin{equation*}
V_{k}:=\left\{ 
\begin{array}{cc}
\limfunc{Sym}^{k}\mathbb{F}^{2} & \text{if }k\geq 0 \\ 
0 & \text{if }k=-1 \\ 
-e^{1+k}\cdot V_{-k-2} & \text{if }k\leq -2%
\end{array}%
.\right.
\end{equation*}
\end{definition}

\bigskip

Notice that the definition of $V_{k}$ for $k<0$ is also suggested by the
expression of the Brauer character of the representation $V_{-k}$ (cf. next
paragraph).

\bigskip

We remark that in the rest of the paper, for any non-negative integer $k$,
the symbol $V_{k}$ will refer sometimes to the $G$-module $\limfunc{Sym}^{k}%
\mathbb{F}^{2}$ and sometimes to its image in $K_{0}\left( G\right) $: the
meaning of the symbol $V_{k}$ will be clear from the context.

\subsection{An identity in $K_{0}\left( G\right) $}

Let us fix an embedding $\iota :\mathbb{F}_{q^{2}}\mathbb{\hookrightarrow }%
M_{2}\left( \mathbb{F}\right) $, corresponding to a choice of $\mathbb{F}$%
-basis for the degree $2$ extension of $\mathbb{F}$ inside $\overline{%
\mathbb{F}}$. Let $\overline{%
\mathbb{Q}
}_{p}$\ be a fixed algebraic closure of the $p$-adic field $%
\mathbb{Q}
_{p}$ and let us fix an isomorphism between $\overline{\mathbb{F}}$ and the
residue field of the ring of integers $\overline{%
\mathbb{Z}
}_{p}$ of $\overline{%
\mathbb{Q}
}_{p}$; denoting by $\chi :\overline{\mathbb{F}}^{\times }\hookrightarrow
\mu \left( \overline{%
\mathbb{Z}
}_{p}\right) $ the corresponding Teichm\"{u}ller lifting, one can compute
the Brauer character $G_{reg}\rightarrow \overline{%
\mathbb{Q}
}_{p}$ of the representations $V_{k}$ ($k\geq 1$) and get:%
\begin{eqnarray*}
\left( 
\begin{array}{cc}
a & 0 \\ 
0 & a%
\end{array}%
\right) &\mapsto &(k+1)\chi \left( a\right) ^{k},\text{ \ \ \ \ \ \ \ \ \ \
\ \ \ \ \ }a\in \mathbb{F}^{\times } \\
\left( 
\begin{array}{cc}
a & 0 \\ 
0 & b%
\end{array}%
\right) &\mapsto &\frac{\chi \left( a\right) ^{k+1}-\chi \left( b\right)
^{k+1}}{\chi \left( a\right) -\chi \left( b\right) },\text{ \ \ \ \ \ }%
a,b\in \mathbb{F}^{\times },a\neq b \\
\iota \left( c\right) &\mapsto &\frac{\chi \left( c\right) ^{q\left(
k+1\right) }-\chi \left( c\right) ^{k+1}}{\chi \left( c\right) ^{q}-\chi
\left( c\right) },\text{ \ \ \ }c\in \mathbb{F}_{q^{2}}^{\times }\backslash 
\mathbb{F}^{\times }.
\end{eqnarray*}

Using the above formulae, in \cite{Serre} Serre deduces:

\begin{proposition}
\label{identit}The following relation holds in $K_{0}\left( G\right) $, for
any integer $k$:%
\begin{equation}
V_{k}-e\cdot V_{k-(q+1)}=V_{k-\left( q-1\right) }-e\cdot V_{k-2q}.
\label{Serre}
\end{equation}
\end{proposition}

\bigskip

More can be shown: let us assume $k>q$ and let $\theta _{q}=X^{q}Y-XY^{q}\in 
\mathbb{F}\left[ X,Y\right] _{q+1}.$

\noindent (Notice this polynomial naturally appears in the classical theory
of the $SL_{2}\left( \mathbb{F}\right) $-invariants of the symmetric algebra 
$\limfunc{Sym}^{\ast }\mathbb{F}^{2}$: by \cite{D}, the ring of $%
SL_{2}\left( \mathbb{F}\right) $-invariants in $\limfunc{Sym}^{\ast }\mathbb{%
F}^{2}$ is generated by the two polynomials%
\begin{equation*}
\theta _{q}=\det \left( 
\begin{array}{cc}
X^{q} & Y^{q} \\ 
X & Y%
\end{array}%
\right) \text{ and }\det \left( 
\begin{array}{cc}
X^{q^{2}} & Y^{q^{2}} \\ 
X & Y%
\end{array}%
\right) /\theta _{q}\text{).}
\end{equation*}%
Let us denote by $\theta _{q}$ also the $G$-equivariant map $e\otimes
V_{k-\left( q+1\right) }\rightarrow V_{k}$ given by multiplication by the
above polynomial. This map is monic and its cokernel is isomorphic to the
induced representation $\limfunc{Ind}_{B}^{G}\left( \eta ^{k}\right) $,
where $B\leq G$ is the subgroup of upper triangular matrices, and $\eta $ is
the character of $B$ given by $\left( 
\begin{array}{cc}
a & b \\ 
0 & c%
\end{array}%
\right) \mapsto a$. Since $\eta $ has order $q-1$ one deduces the following
result (cf. \ref{Serre}):

\begin{proposition}
\label{sonno}For $k>q$, there is an exact sequence of $G$-modules $%
0\rightarrow e\otimes V_{k-\left( q+1\right) }\rightarrow V_{k}\rightarrow 
\limfunc{Ind}_{B}^{G}\left( \eta ^{k}\right) \rightarrow 0$; in particular,
for $k\geq 2q$ we have a $G$-isomorphism:%
\begin{equation}
\frac{V_{k}}{e\otimes V_{k-(q+1)}}\simeq \frac{V_{k-\left( q-1\right) }}{%
e\otimes V_{k-2q}}  \label{ciaoo}
\end{equation}
\end{proposition}

\bigskip

We can obtain the isomorphism \ref{ciaoo} also as follows: as in \cite{A-S},
we let $I_{k}$ be,\ for any $k\geq 0$, the space of functions $\mathbb{F}%
^{2}\rightarrow \mathbb{F}$ that are homogeneous of degree $k$ and vanish at 
$(0,0)$, endowed with the $G$-action given by $\left( g\cdot f\right) \left(
x,y\right) =f((g^{-1}(x,y)^{t})^{t})$ for $g\in G,f\in I_{k}$ and $(x,y)\in 
\mathbb{F}^{2}$. If $\rho :G\rightarrow GL\left( V\right) $ is any
representation of $G$ over some field, we denote by $\rho ^{t}$ the
representation defined by $\rho ^{t}\left( g\right) :=\rho (\left(
g^{t}\right) ^{-1})$, where $g^{t}$ is the transpose of the matrix $g\in G$.
Denote then by $\tau _{k}:V_{k}\rightarrow I_{k}^{t}$ the $G$-map sending a
polynomial in $V_{k}$ to the associated polynomial function in $I_{k}^{t}$.

Assume now $k>q$ and let $T_{O}:=X^{k-\left( q-1\right)
}\tprod\nolimits_{a\in \mathbb{F}^{\times }}\left( X-aY\right) \in V_{k}$:
then $\tau _{k}\left( T_{O}\right) $ is non-zero only at $\left[ 1:0\right]
\in \mathbb{P}^{1}\left( \mathbb{F}\right) $; since $G$ act transitively on $%
\mathbb{P}^{1}\left( \mathbb{F}\right) $, for any $P\in \mathbb{P}^{1}\left( 
\mathbb{F}\right) $ we can find a polynomial $T_{P}\in V_{k}$ such that $%
\tau _{k}\left( T_{P}\right) $ is non zero only at $P$:\ $\tau _{k}$ is
therefore onto. Since $\theta _{q}\left( e\otimes V_{k-\left( q+1\right)
}\right) \leq \ker \tau _{k}$ and $\dim _{\mathbb{F}}I_{k}^{t}=q+1$ we
conclude that $\limfunc{coker}\left( \theta _{q}:e\otimes V_{k-\left(
q+1\right) }\rightarrow V_{k}\right) \simeq I_{k}^{t}=I_{k+(q-1)}^{t}.$

Notice that for any $k\geq 0$ we have $I_{k}\simeq e^{-k}\otimes \mathbb{F}%
\left[ \mathbb{P}^{1}\left( \mathbb{F}\right) \right] $, where $\mathbb{F}%
\left[ \mathbb{P}^{1}\left( \mathbb{F}\right) \right] $ is the $G$-module of 
$\mathbb{F}$-valued functions on $\mathbb{P}^{1}\left( \mathbb{F}\right) $
with the usual left $G$-action; denoting by $M_{0}$ (resp $M_{1}$) the $G$%
-submodule of $\mathbb{F}\left[ \mathbb{P}^{1}\left( \mathbb{F}\right) %
\right] $ consisting of the constant (resp. average zero) functions, we have 
$\mathbb{F}\left[ \mathbb{P}^{1}\left( \mathbb{F}\right) \right]
=M_{0}\oplus M_{1}.$ We now need two lemmas:

\begin{lemma}
In $\mathbb{F}\left[ X,Y\right] $ one has the identity $X^{q-1}+\sum%
\nolimits_{a\in \mathbb{F}}\left( aX+Y\right) ^{q-1}=0.$
\end{lemma}

\begin{proof}
Writing $\tsum\nolimits_{a\in \mathbb{F}}\left( aX+Y\right)
^{q-1}=\tsum\nolimits_{a\in \mathbb{F}}\sum_{j=0}^{q-1}\tbinom{q-1}{j}%
a^{j}X^{j}Y^{q-1-j}$, we see that the coefficient of $X^{q-1}$ (resp. $%
Y^{q-1}$) in this polynomial is $q-1$ (resp. $0$), hence it's enough we
prove that for any $0<j<q-1$ we have $\sum_{a\in \mathbb{F}^{\times
}}a^{j}=0.$ Let us chose a generator $\varepsilon $ of the multiplicative
group of $\mathbb{F}$ and let $\lambda $ be the order of $\varepsilon ^{j}$
in this group; if $q-1=\lambda d$ in $%
\mathbb{Z}
$, we have:%
\begin{equation*}
\tsum\nolimits_{a\in \mathbb{F}^{\times
}}a^{j}=\tsum\nolimits_{n=0}^{q-2}\varepsilon
^{nj}=d\tsum\nolimits_{n=0}^{\lambda -1}\varepsilon ^{nj}.
\end{equation*}

\noindent Since $0=\varepsilon ^{\lambda j}-1=\left( \varepsilon
^{j}-1\right) \sum_{n=0}^{\lambda -1}\varepsilon ^{nj}$ and $\varepsilon
^{j}\neq 1$, we conclude that the above summation is zero.
\end{proof}

\begin{lemma}
\label{questooo}There exists a $G$-module epimorphism $\vartheta :\mathbb{F}%
\left[ \mathbb{P}^{1}\left( \mathbb{F}\right) \right] \rightarrow V_{q-1}$
whose kernel is $M_{0}$. In particular $M_{1}\simeq V_{q-1}$ as $G$-modules.
\end{lemma}

\begin{proof}
Let $\vartheta $ be the $G$-map defined as follows:%
\begin{equation*}
f\in \mathbb{F}\left[ \mathbb{P}^{1}\left( \mathbb{F}\right) \right] \text{
\ \ }\mapsto \text{ \ }\tsum\nolimits_{\left[ a:b\right] \in \mathbb{P}%
^{1}\left( \mathbb{F}\right) }f\left( \left[ a:b\right] \right) \left(
aX+bY\right) ^{q-1}.
\end{equation*}

For any point $P$ of the $\mathbb{F}$-projective line, let $f_{P}\in \mathbb{%
F}\left[ \mathbb{P}^{1}\left( \mathbb{F}\right) \right] $ be the function
equal to $1$ at $P$ and zero everywhere else. Notice that $X^{q-1}=\vartheta
\left( f_{\left[ 1:0\right] }\right) $; let then $\alpha $ be an integer, $%
0\leq \alpha <q-1$ and let $f=\sum_{P\in \mathbb{P}^{1}\left( \mathbb{F}%
\right) }\lambda _{P}f_{P}$, where $\lambda _{P}\in \mathbb{F}$. We have $%
\vartheta \left( f\right) =X^{\alpha }Y^{q-1-\alpha }$ if, for example, $%
\lambda _{\left[ 1:0\right] }=0$ and the following $q$\ equations in the
variables $\left\{ \lambda _{\left[ a:1\right] }\right\} _{a\in \mathbb{F}}$
have a common solution:%
\begin{eqnarray*}
\tsum\nolimits_{a\in \mathbb{F}}\lambda _{\left[ a:1\right] }a^{\alpha } &=&%
\tbinom{q-1}{\alpha }^{-1} \\
\tsum\nolimits_{a\in \mathbb{F}}\lambda _{\left[ a:1\right] }a^{j} &=&0\text{
\ \ \ \ \ (}0\leq j\leq q-1,j\neq \alpha \text{).}
\end{eqnarray*}

\noindent (Notice that the identity $(X+Y)^{q}=\left( X+Y\right)
\sum_{i=0}^{q-1}\left( -1\right) ^{i}X^{q-1-i}Y^{i}$ in $\mathbb{F}\left[ X,Y%
\right] $\ implies that $\tbinom{q-1}{i}\equiv (-1)^{i}(\func{mod}p)$ for
any $0\leq i\leq q-1$, hence $\tbinom{q-1}{\alpha }$ is invertible in $%
\mathbb{F}$).

Since the matrix of this $q\times q$ system is of Vandermonde type over the
elements of $\mathbb{F}$, the system has a unique solution, so that $%
\vartheta ^{-1}\left( X^{\alpha }Y^{q-1-\alpha }\right) \neq \varnothing $
and $\vartheta $ is onto. Since the kernel of $\vartheta $ has dimension $1$%
, the previous lemma implies $\ker \vartheta =M_{0}$.
\end{proof}

As a consequence of our computations we conclude:

\begin{proposition}
\label{fame}\noindent For any integer $k>q$ and any $\lambda \geq 0$\ there
are isomorphisms of $G$-modules:%
\begin{equation*}
\frac{V_{k}}{e\otimes V_{k-(q+1)}}\simeq e^{k}\otimes \left( \mathbb{F\oplus 
}V_{q-1}^{t}\right) \simeq \frac{V_{k+\lambda (q-1)}}{e\otimes V_{k+\lambda
(q-1)-(q+1)}},
\end{equation*}

\noindent where the inclusion $e\otimes V_{k+\lambda
(q-1)-(q+1)}\hookrightarrow V_{k+\lambda (q-1)}$ is induced by the
multiplication by $\theta _{q}$.
\end{proposition}

\begin{proof}
Just notice that $I_{k}^{t}\simeq \left( e^{-k}\otimes \left( \mathbb{%
F\oplus }V_{q-1}\right) \right) ^{t}\simeq e^{k}\otimes \left( \mathbb{F}^{t}%
\mathbb{\oplus }V_{q-1}^{t}\right) =e^{k}\otimes \left( \mathbb{F\oplus }%
V_{q-1}^{t}\right) $.
\end{proof}

\section{Reduction of cuspidal representations}

Let us rewrite the identity \ref{Serre} as: 
\begin{equation}
V_{k}-V_{k-\left( q-1\right) }=e\cdot \left( V_{k-(q+1)}-V_{k-2q}\right) ,%
\text{ \ \ \ (}k\in 
\mathbb{Z}
\text{).}  \label{questa!!!}
\end{equation}

\noindent \noindent One is naturally led to ask when $V_{k}-V_{k-\left(
q-1\right) }$ is positive in $K_{0}(GL_{2}\left( \mathbb{F}\right) )$: we
would like to have, at least for some values of $k$, a monic map of $%
GL_{2}\left( \mathbb{F}\right) $-modules $V_{k-\left( q-1\right)
}\rightarrow V_{k}$ that "raises the weight by $q-1$". The Serre $D$-map
will play this role, if $1\leq k\leq p-1$.

\subsection{The $D$-map of Serre}

We assume from now on that $k\geq 0$ and we will impose some other
restrictions on it later.

\noindent Let $\chi :\mathbb{F}_{q^{2}}^{\times }\rightarrow \overline{%
\mathbb{Q}
}_{p}^{\times }$ be the (restriction of the) Teichm\"{u}ller character fixed
in the previous section, and let $k$ be an integer such that $k\neq 0(\func{%
mod}q+1)$, so that $\chi ^{k}$ is indecomposable (i.e. it does not factor
through the norm map $\mathbb{F}_{q^{2}}^{\times }\rightarrow \mathbb{F}%
^{\times }$). Under these assumptions, there is a unique $\overline{%
\mathbb{Q}
}_{p}$-representation $\Theta \left( \chi ^{k}\right) $ of $G$ characterized
by the property that $\Theta \left( \chi ^{k}\right) \otimes \limfunc{sp}%
\simeq \limfunc{Ind}_{\mathbb{F}_{q^{2}}^{\times }}^{G}\left( \chi
^{k}\right) $, where $\mathbb{F}_{q^{2}}^{\times }$ is embedded in $G$ by
the map $\iota $ fixed in the previous section, and $\limfunc{sp}$ is the $G$%
-subrepresentation of $\overline{%
\mathbb{Q}
}_{p}\left[ \mathbb{P}^{1}\left( \mathbb{F}\right) \right] $ consisting of
functions of average zero. $\Theta \left( \chi ^{k}\right) $ is the cuspidal
representation associated to $\chi ^{k}$.

\noindent The character of $\Theta \left( \chi ^{k}\right) $ is given by:

\begin{eqnarray*}
\left( 
\begin{array}{cc}
a & 0 \\ 
0 & a%
\end{array}%
\right) &\mapsto &(q-1)\chi \left( a\right) ^{k},\text{ \ \ \ \ \ \ \ \ \ }%
a\in \mathbb{F}^{\times } \\
\left( 
\begin{array}{cc}
a & 1 \\ 
0 & a%
\end{array}%
\right) &\mapsto &-\chi \left( a\right) ^{k},\text{ \ \ \ \ \ \ \ \ \ \ \ \
\ \ \ \ \ \ }a\in \mathbb{F}^{\times } \\
\left( 
\begin{array}{cc}
a & 0 \\ 
0 & b%
\end{array}%
\right) &\mapsto &0,\text{ \ \ \ \ \ \ \ \ \ \ \ \ \ \ \ \ \ \ \ }a,b\in 
\mathbb{F}^{\times },a\neq b \\
\iota \left( c\right) &\mapsto &-\chi \left( c\right) ^{k}-\chi \left(
c^{q}\right) ^{k},\text{ \ \ \ \ }c\in \mathbb{F}_{q^{2}}^{\times
}\backslash \mathbb{F}^{\times }.
\end{eqnarray*}

A computation of Brauer characters gives:

\begin{proposition}
Let $k$ be an integer such that $k\neq 0(\func{mod}q+1)$; the following
identity holds in $K_{0}\left( G\right) :$ 
\begin{equation*}
V_{k+\left( q-1\right) }-V_{k}=\overline{\widetilde{\Theta }\left( \chi
^{k}\right) },
\end{equation*}%
\noindent where the bar on the left hand side denotes the reduction $\func{%
mod}p$ of \textit{any} fixed integral model $\widetilde{\Theta }\left( \chi
^{k}\right) $ of $\Theta \left( \chi ^{k}\right) $ (the choice of integral
model is not relevant since the identity is written in $K_{0}\left( G\right) 
$).
\end{proposition}

\bigskip

\begin{definition}
\label{polar}\textbf{(J-P. Serre)} The map\footnote{%
An analogue of the map $D$ in characteristic zero appears in the theory of $%
G $-invariants, in particular cf. the definition of the polarization map in 
\cite{HaaJan}, Appendix F, \S 1.} $D:\mathbb{F}\left[ X,Y\right] \rightarrow 
\mathbb{F}\left[ X,Y\right] $ is the $\mathbb{F}$-linear map defined by:%
\begin{equation*}
Df\left( X,Y\right) :=X^{q}\cdot \frac{\partial f}{\partial X}+Y^{q}\cdot 
\frac{\partial f}{\partial Y},
\end{equation*}

\noindent for any polynomial $f\in \mathbb{F}\left[ X,Y\right] $ (notice
here $q=\#\mathbb{F}$).
\end{definition}

The map $D$ is $G$-equivariant: if $%
\begin{pmatrix}
a & b \\ 
c & d%
\end{pmatrix}%
\in G$ and $\alpha ,\beta \geq 0$ the polynomial $D\left( g\cdot X^{\alpha
}Y^{\beta }\right) $ equals:%
\begin{equation*}
\alpha \left( aX+cY\right) ^{q}\left( aX+cY\right) ^{\alpha -1}\left(
bX+dY\right) ^{\beta }+\beta \left( bX+dY\right) ^{q}\left( aX+cY\right)
^{\alpha }\left( bX+dY\right) ^{\beta -1}.
\end{equation*}%
One easily sees that this is exactly $\left( gX\right) ^{q}\left( g\cdot
\alpha X^{\alpha -1}Y^{\beta }\right) +\left( gY\right) ^{q}\left( g\cdot
\beta X^{\alpha }Y^{\beta -1}\right) $.

\noindent Furthermore, N. Fakhruddin proved the following:

\begin{proposition}
\label{injective}The kernel of the map $D:\mathbb{F}\left[ X,Y\right]
\rightarrow \mathbb{F}\left[ X,Y\right] $ is given by $\ker D=\mathbb{F}%
\left[ X^{p},Y^{p},\theta _{q}\right] .$
\end{proposition}

\begin{proof}
Let $A=\mathbb{F[}X^{p},Y^{p},\theta _{q}\mathbb{]}$ and $B=\ker D$; notice
that $B$ is a ring and we have the inclusions $\mathbb{F}[X^{p},Y^{p}]%
\subseteq A\subseteq B\subseteq \mathbb{F}[X,Y]$. The polynomial $%
t^{p}-(X^{pq}Y^{p}-X^{p}Y^{pq})\in \mathbb{F}(X^{p},Y^{p})[t]$ is
irreducible in $\mathbb{F}(X^{p},Y^{p})[t]$ since $X^{pq}Y^{p}-X^{p}Y^{pq}$
does not have a $p^{th}$-root in $\mathbb{F}(X^{p},Y^{p})$, so that we have $%
[Q\left( A\right) :\mathbb{F}(X^{p},Y^{p})]=p$, where we denote by $Q(R)$
the field of fractions of an integral domain $R$ inside some extension of $R$%
.

Now observe that $Q\left( B\right) $ is properly contained inside $\mathbb{F}%
(X,Y)$: if not, we could write $X=\frac{f}{g}$ with $f,g\in B$, $g\neq 0$
and $1=af+bg$ for some $a,b\in B$; this would imply $X=f\cdot \left(
aX+b\right) $ so that $f$ $\in B$\ would be an associate of $X$ in $\mathbb{F%
}[X,Y]$. Since $[\mathbb{F}(X,Y):\mathbb{F}(X^{p},Y^{p})]=p^{2}$ we have
therefore $Q\left( A\right) =Q\left( B\right) $. Notice that $\mathbb{F}%
[X,Y]/A$ is an integral extension, so that $B/A$ is too.

Now observe that the domain $A$ is normal, since the corresponding variety
has equation $X_{1}^{q}X_{2}-X_{1}X_{2}^{q}-X_{3}^{p}=0$, and then it
defines an hypersurface of $\mathbb{A}_{/\mathbb{F}}^{3}$ that is
non-singular in codimension one (cf. \cite{Har}, II.8.23). We conclude $A=B$.
\end{proof}

\bigskip

We will need the following:

\begin{lemma}
If $q=p>3$ is a prime number, then $V_{p}$ is a non-trivial extension of $%
V_{p-2}$ via $V_{1}$ in the category of $\mathbb{F}_{p}\left[ SL_{2}(\mathbb{%
F}_{p})\right] $-modules. In particular, for $p>3$,\ the one dimensional $%
\mathbb{F}_{p}$-space $\limfunc{Ext}_{SL_{2}\left( \mathbb{F}_{p}\right)
}^{1}\left( V_{p-2},V_{1}\right) \simeq \limfunc{Ext}_{SL_{2}\left( \mathbb{F%
}_{p}\right) }^{1}\left( V_{1},V_{p-2}\right) $ is generated by the class of 
$V_{p}$.
\end{lemma}

\begin{proof}
A computation of Brauer characters tells us that $V_{p}-V_{1}=V_{p-2}$ in $%
K_{0}\left( SL_{2}(\mathbb{F}_{p})\right) $; the existence of the monic map $%
D:V_{1}\hookrightarrow V_{p}$ and the fact that $V_{p-2}$ is irreducible
give therefore an exact sequence of the form:%
\begin{equation}
0\rightarrow V_{1}\overset{D}{\rightarrow }V_{p}\overset{\pi ^{\prime }}{%
\rightarrow }V_{p-2}\rightarrow 0.  \label{uniii}
\end{equation}

We can choose $\pi ^{\prime }$ so that it sends $X^{p}$ and $Y^{p}$ to zero
and $\pi ^{\prime }(X^{r}Y^{p-r}):=k\left( r\right) \cdot X^{r-1}Y^{p-1-r}$
for $0<r<p$, where $k(r):=\frac{\tbinom{p-2}{r-1}}{\tbinom{p-1}{r}}\neq 0$
(cf. the proof of Proposition \ref{uffaaa}). With this choice of $\pi
^{\prime }$, the sequence:%
\begin{equation}
0\rightarrow V_{1}\overset{D}{\rightarrow }V_{p}\overset{\pi ^{\prime }}{%
\rightarrow }e\otimes V_{p-2}\rightarrow 0  \label{duiiii}
\end{equation}%
is exact in the category of $GL_{2}(\mathbb{F}_{p})$-modules; since $SL_{2}(%
\mathbb{F}_{p})$ is normal in $GL_{2}(\mathbb{F}_{p})$ of index prime to $p$%
, we have that \ref{uniii} splits as a sequence of $SL_{2}(\mathbb{F}_{p})$%
-modules if and only if \ref{duiiii} splits as a sequence of $GL_{2}(\mathbb{%
F}_{p})$-modules. We will show that \ref{duiiii}\ does not split.

A splitting $GL_{2}(\mathbb{F}_{p})$-homomorphism $\sigma :e\otimes
V_{p-2}\rightarrow V_{p}$ for $\pi ^{\prime }$\ has to send $k\left(
r\right) \cdot X^{r-1}Y^{p-1-r}$ to $X^{r}Y^{p-r}+h_{1}\left( r\right)
X^{p}+h_{2}\left( r\right) Y^{p}$ for some $h_{1}\left( r\right)
,h_{2}\left( r\right) \in \mathbb{F}_{p}$ ($0<r<p$); imposing the
equivariance of the map $\sigma $ with respect to the maximal split torus of
diagonal matrices in $GL_{2}(\mathbb{F}_{p})$, we find the conditions: 
\begin{eqnarray*}
h_{1}\left( r\right) &=&a^{r-1}b^{1-r}\cdot h_{1}(r), \\
h_{2}(r) &=&a^{r}b^{-r}\cdot h_{2}\left( r\right) ,
\end{eqnarray*}
that have to be satisfied for every $a,b\in \mathbb{F}_{p}^{\times }$ and
every $0<r<p$. We therefore have $h_{1}(r)=0$ if $2\leq r\leq p-1$ and $%
h_{2}(r)=0$ if $1\leq r\leq p-2$.

This implies a contradiction, since then - being $p>3$ - we would get:%
\begin{equation*}
\sigma \left( 
\begin{pmatrix}
1 & 0 \\ 
1 & 1%
\end{pmatrix}%
\cdot k(2)XY^{p-3}\right) \neq 
\begin{pmatrix}
1 & 0 \\ 
1 & 1%
\end{pmatrix}%
\cdot \sigma \left( k(2)XY^{p-3}\right) .
\end{equation*}

The second assertion in the Lemma follows from the Propositions in \S \S %
12.1,12.2 of \cite{Hum}.
\end{proof}

\bigskip

We are now able to prove the following result, due to G. Savin:

\begin{proposition}
\label{essential}If $q=p>3$, then there are no $\mathbb{F}_{p}\left[ GL_{2}(%
\mathbb{F}_{p})\right] $-embeddings $V_{p}\rightarrow V_{2p-1}.$
\end{proposition}

\begin{proof}
The $\mathbb{F}_{p^{2}}[SL_{2}(\mathbb{F}_{p^{2}})]$-module $\left(
V_{p-1}\otimes _{\mathbb{F}_{p}}\mathbb{F}_{p^{2}}\right) \otimes _{\mathbb{F%
}_{p^{2}}}\left( V_{1}\otimes _{\mathbb{F}_{p}}\mathbb{F}_{p^{2}}\right)
^{Fr}$ is irreducible by Steinberg's tensor product theorem, and it is
isomorphic to $V_{2p-1}\otimes _{\mathbb{F}_{p}}\mathbb{F}_{p^{2}}$ via the
map induced by $X^{i}Y^{p-1-i}\otimes X^{j}Y^{1-j}\mapsto
X^{i+pj}Y^{2p-1-i-pj}$ ($0\leq i\leq p-1$, $0\leq j\leq 1$). We deduce the
existence of an isomorphism of $\mathbb{F}_{p}\left[ SL_{2}(\mathbb{F}_{p})%
\right] $-modules $V_{2p-1}\simeq V_{1}\otimes _{\mathbb{F}_{p}}V_{p-1}$.

Since $V_{1}$ is an $\mathbb{F}_{p}\left[ SL_{2}(\mathbb{F}_{p})\right] $%
-submodule of $V_{p}$, to prove the proposition it suffices to show:%
\begin{equation*}
\limfunc{Hom}\nolimits_{\mathbb{F}_{p}[SL_{2}(\mathbb{F}%
_{p})]}(V_{1},V_{2p-1})=0.
\end{equation*}
We have isomorphisms of $\mathbb{F}_{p}$-vector spaces:%
\begin{eqnarray*}
\limfunc{Hom}\nolimits_{\mathbb{F}_{p}[SL_{2}(\mathbb{F}%
_{p})]}(V_{1},V_{2p-1}) &\simeq &\limfunc{Hom}\nolimits_{\mathbb{F}%
_{p}[SL_{2}(\mathbb{F}_{p})]}(V_{1},V_{1}\otimes _{\mathbb{F}_{p}}V_{p-1}) \\
&\simeq &\limfunc{Hom}\nolimits_{\mathbb{F}_{p}[SL_{2}(\mathbb{F}%
_{p})]}(V_{1}\otimes _{\mathbb{F}_{p}}V_{1},V_{p-1}),
\end{eqnarray*}

\noindent where the second isomorphism is a consequence of the self-duality
of $V_{1}$ as $\mathbb{F}_{p}[SL_{2}(\mathbb{F}_{p})]$-module. Now notice
that the vector space $\limfunc{Hom}\nolimits_{\mathbb{F}_{p}[SL_{2}(\mathbb{%
F}_{p})]}(V_{1}\otimes _{\mathbb{F}_{p}}V_{1},V_{p-1})$ is trivial, since $%
V_{p-1}$ is a simple $\mathbb{F}_{p}\left[ SL_{2}(\mathbb{F}_{p})\right] $%
-module of dimension $p>3$.
\end{proof}

\bigskip

The last proposition shows that when $q=p>3$ and the Serre $D$-map fails to
be injective as map $V_{p}\rightarrow V_{2p-1}$, then there are actually no
subobjects of $V_{2p-1}$ isomorphic to $V_{p}$. Notice that an analogous
result is clearly true for maps $V_{0}\rightarrow V_{p-1}$, since $%
V_{p-1}=St $ is irreducible over $\mathbb{F}_{p}\mathbb{[}GL_{2}\left( 
\mathbb{F}_{p}\right) \mathbb{]}$.

\bigskip

We also have:

\begin{proposition}
\label{negativity}The element $V_{k}-V_{k-(q-1)}$ of $K_{0}(G)$ is positive
if and only if $k\neq -2(\func{mod}q+1)$.
\end{proposition}

\begin{proof}
By \ \ref{questa!!!}\ we can restrict ourselves to the case $-2\leq k\leq
q-2 $. If $0\leq k\leq q-2$ then $k-(q-1)<0$, so that $%
V_{k}-V_{k-(q-1)}=V_{k}-V_{-(q-1-k)}=V_{k}-\left( -e^{1-(q-1-k)}\cdot
V_{(q-1-k)-2}\right) =V_{k}+e^{1+k}\cdot V_{q-k-3}$, and this element is
positive in $K_{0}(G)$ (notice that here both $k$ and $q-k-3$ are
non-negative integers). If $k=-1$, then $%
V_{-1}-V_{-1-(q-1)}=0+V_{q-2}=V_{q-2}$ is again positive in the Grothendieck
group of $G$.

If $k=-2$, we need to consider the element $V_{-2}-V_{-(q+1)}$ of $K_{0}(G)$%
: this is equal to $e^{-1}\cdot (V_{q-1}-V_{0})$ and is not positive, since $%
V_{q-1}=St$ and $V_{0}$ are distinct non-zero irreducible representations of 
$G$.
\end{proof}

\bigskip

\bigskip

Assuming now that $1\leq k\leq p-1$, the map of $G$-modules $%
D:V_{k}\rightarrow V_{k+(q-1)}$ is monic and we deduce, via a computation of
Brauer characters, the existence of a $G$-module isomorphism:

\begin{center}
\begin{equation}
\left( \dfrac{V_{k+(q-1)}}{D\left( V_{k}\right) }\right) ^{ss}\simeq \left( 
\overline{\widetilde{\Theta }\left( \chi ^{k}\right) }\right) ^{ss},
\label{key!}
\end{equation}
\end{center}

\noindent where for a $G$-module $M$, $M^{ss}$ denotes the $G$-module
obtained by semi-simplifying $M$.

We will show that this isomorphism comes from an isomorphism of modules
without the semi-simplification, for $2\leq k\leq p-1$, $k\neq \frac{q+1}{2}$
and some integral model $\widetilde{\Theta }\left( \chi ^{k}\right) $ of $%
\Theta \left( \chi ^{k}\right) $.

\subsection{An exact sequence from crystalline cohomology}

We recall some constructions from \cite{HaaJan}.\ Let $R$ be either a field
extension $\mathbb{K}$ of $\mathbb{F}_{q^{2}}$ inside $\overline{\mathbb{F}}$%
\ or the ring of Witt vectors $W=W\left( \mathbb{K}\right) $ of such a
field, together with the natural embedding $\mathbb{F}_{q^{2}}^{\times
}\hookrightarrow R^{\times }$; let $A=\left( 
\begin{array}{cc}
0 & 1 \\ 
-1 & 0%
\end{array}%
\right) \in SL_{2}\left( R\right) $ and:

\begin{equation*}
H\left( X,Y,Z\right) =\left( X,Y,Z\right) \cdot \left( 
\begin{array}{cc}
A & 0 \\ 
0 & -1%
\end{array}%
\right) \cdot \overline{\left( X,Y,Z\right) }^{t}\in R\left[ X,Y,Z\right] ,
\end{equation*}%
where the upper bar denotes the $q^{th}$-power map. The homogeneous
polynomial $H$ is irreducible over $R$ and\ defines a projective, smooth and
irreducible curve $C_{/R}$ having equation $XY^{q}-X^{q}Y-Z^{q+1}=0$.

We let $U_{2}\left( \mathbb{F}_{q^{2}}\right) =\left\{ g\in GL_{2}\left( 
\mathbb{F}_{q^{2}}\right) :g\cdot A\cdot \overline{g}^{t}=A\right\} $ and we
consider $\mathbb{F}_{q^{2}}^{\times }$ embedded in $U_{2}\left( \mathbb{F}%
_{q^{2}}\right) $ via $t\mapsto \left( 
\begin{array}{cc}
t & 0 \\ 
0 & t^{-q}%
\end{array}%
\right) $, so that $U_{2}\left( \mathbb{F}_{q^{2}}\right) =\mathbb{F}%
_{q^{2}}^{\times }\cdot SL_{2}\left( \mathbb{F}\right) $. By construction, $%
U_{2}\left( \mathbb{F}_{q^{2}}\right) $ acts on $C_{/\mathbb{K}}$ via the
embedding:

\begin{equation*}
U_{2}\left( \mathbb{F}_{q^{2}}\right) \hookrightarrow GL_{3}\left( \mathbb{F}%
_{q^{2}}\right) :g\mapsto \left( 
\begin{array}{cc}
g & 0 \\ 
0 & 1%
\end{array}%
\right) ,
\end{equation*}%
and the induced action of $\mathbb{F}_{q^{2}}^{\times }\hookrightarrow
U_{2}\left( \mathbb{F}_{q^{2}}\right) $ lifts to an action on $C_{/W}$. In
particular the group $\mu =\left\{ t\in \mathbb{F}_{q^{2}}^{\times
}:t^{q+1}=1\right\} \subset \mathbb{F}_{q^{2}}^{\times }$ acts on $C_{/R}$
and its cohomology groups.

\bigskip

Let $i\geq 0$ be any integer; we denote by $H_{cris}^{i}\left( C_{/\mathbb{K}%
}\right) $ the $i^{th}$ group of the crystalline cohomology of $C_{/\mathbb{K%
}}$ (cf. \cite{Ill}), and we let $H_{dR}^{i}\left( C_{/R}\right) $ be the $%
i^{th}$ group of the de Rham cohomology of $C_{/R}$. Since $C_{/\mathbb{K}}$
admits a smooth lifting to characteristic zero (i.e. $C_{/W}$), we have
canonical isomorphisms of $W$-modules $H_{cris}^{i}\left( C_{/\mathbb{K}%
}\right) \simeq H_{dR}^{i}\left( C_{/W}\right) $.

The de Rham cohomology groups of $C_{/R}$ are computed as hypercohomology of
the de Rham complex for $C_{/R}$, so that we have the corresponding
(Hodge-de Rham) spectral sequence $E_{1}^{p,q}=H^{q}(C_{/R},\Omega
_{C_{/R}}^{p})\Longrightarrow H_{dR}^{p+q}(C_{/R})$ ($p,q\geq 0$) that
degenerates at the first page $E_{1}$, since $C_{/R}$ is a curve (cf. \cite%
{Dn}, II.6.). We deduce from this that $H_{dR}^{i}\left( C_{/R}\right)
\simeq R$ for $i=0,2$ and $H_{dR}^{i}\left( C_{/R}\right) =0$ for $i>2$.

For $i=1$ the degeneracy of $E_{1}^{p,q}$ at the first page gives:

\begin{proposition}
\label{here}There is a natural exact sequence of free $R$-modules:%
\begin{equation*}
0\rightarrow H^{0}(C_{/R},\Omega _{C_{/R}}^{1})\rightarrow
H_{dR}^{1}(C_{/R})\rightarrow H^{1}(C_{/R},\mathcal{O}_{C_{/R}})\rightarrow
0.
\end{equation*}
\end{proposition}

For $1\leq k\leq q$ let $k:\mu \rightarrow \mu $ denote the character $%
t\mapsto t^{k}$ ($t\in \mu $). Since $\left\vert \mu \right\vert =q+1$ is
invertible in $R$, the proposition gives:

\begin{corollary}
For every $1\leq k\leq q$, there is a natural exact sequence of $R$-modules,
:%
\begin{equation*}
0\rightarrow H^{0}(C_{/R},\Omega _{C_{/R}}^{1})_{-k}\rightarrow
H_{dR}^{1}(C_{/R})_{-k}\rightarrow H^{1}(C_{/R},\mathcal{O}%
_{C_{/R}})_{-k}\rightarrow 0,
\end{equation*}

\noindent where the subscript $-k$ designates the corresponding $\mu $%
-eigenspace.
\end{corollary}

In the particular case $R=\mathbb{K}$, we have an action of $U_{2}\left( 
\mathbb{F}_{q^{2}}\right) $ on $C_{/\mathbb{K}}$ and its cohomology. We have:

\begin{proposition}
\label{ES1}There is a natural exact sequence of $\mathbb{K}\left[
U_{2}\left( \mathbb{F}_{q^{2}}\right) \right] $-modules, for every $1\leq
k\leq q$:%
\begin{equation*}
0\rightarrow H^{0}(C_{/\mathbb{K}},\Omega _{C_{/\mathbb{K}%
}}^{1})_{-k}\rightarrow H_{dR}^{1}(C_{/\mathbb{K}})_{-k}\rightarrow
H^{1}(C_{/\mathbb{K}},\mathcal{O}_{C_{/\mathbb{K}}})_{-k}\rightarrow 0.
\end{equation*}
\end{proposition}

We will interested in the sequel in writing down explicitly the maps
occurring in the above sequences.

\begin{remark}
\label{reduction}For any $i\geq 0$ there is an exact sequence of "universal
coefficients":%
\begin{equation*}
0\rightarrow H_{cris}^{i}(C_{/\mathbb{K}})\otimes _{W}\mathbb{K\rightarrow }%
H_{dR}^{i}(C_{/\mathbb{K}})\rightarrow \limfunc{Tor}%
\nolimits_{1}^{W}(H_{cris}^{i+1}(C_{/\mathbb{K}}),\mathbb{K)\rightarrow }0.
\end{equation*}

\noindent (cf. \cite{Ill}). Taking $i=1$ and using the identification $%
H_{cris}^{1}\left( C_{/\mathbb{K}}\right) \simeq H_{dR}^{1}\left(
C_{/W}\right) $\ and the fact that $\limfunc{Tor}\nolimits_{1}^{W}(W,\mathbb{%
K)=}0$, we deduce the existence of a natural isomorphism: 
\begin{equation*}
H_{dR}^{1}(C_{/\mathbb{K}})\simeq H_{dR}^{1}\left( C_{/W}\right) \otimes _{W}%
\mathbb{K}.
\end{equation*}
\end{remark}

\subsection{An exact sequence of representations of $U_{2}\left( \mathbb{F}%
_{q^{2}}\right) $}

We work over the field $\mathbb{F}_{q^{2}}$, the quadratic extension of $%
\mathbb{F}$ inside $\overline{\mathbb{F}}$. We assume from now on that $p$
is an \textit{odd} prime and we put as usual $G=GL_{2}\left( \mathbb{F}%
\right) $; let $k$ be an integer such that $2\leq k\leq p-1$, and let $V_{k}=%
\limfunc{Sym}^{k}\mathbb{F}_{q^{2}}^{2}$ be the usual left representation of 
$U_{2}\left( \mathbb{F}_{q^{2}}\right) \subseteq GL_{2}\left( \mathbb{F}%
_{q^{2}}\right) $ over $\mathbb{F}_{q^{2}}$. Notice that we changed the
notation from section $2$: in this section we will be denoting by $%
V_{k}^{^{\prime }}$ the $G$-representation $\limfunc{Sym}^{k}\mathbb{F}^{2}$%
, so that $V_{k}=V_{k}^{^{\prime }}\otimes _{\mathbb{F}}\mathbb{F}%
_{q^{2}}^{2}$, as $G$-representations.

The multiplication by $\theta _{q}$ (not by $\theta _{q^{2}}$!) induces a $%
\mathbb{F}_{q^{2}}\left[ U_{2}\left( \mathbb{F}_{q^{2}}\right) \right] $%
-monomorphism $\theta _{q}:V_{k-2}\hookrightarrow V_{k+q-1}$. Let $D^{\prime
}:V_{k}^{^{\prime }}\hookrightarrow V_{k+(q-1)}^{^{\prime }}$ be the Serre $%
D $-map that increases the degree of $q-1$, as in Definition \ref{polar};
the monic $\mathbb{F}_{q^{2}}$-linear map $D:=D^{^{\prime }}\otimes \mathbb{F%
}_{q^{2}}:V_{k}\hookrightarrow V_{k+(q-1)}$ is $SL_{2}\left( \mathbb{F}%
\right) $-equivariant, but not $\mathbb{F}_{q^{2}}^{\times }$-equivariant;
nevertheless $D\left( V_{k}\right) $ is an $\mathbb{F}_{q^{2}}\left[
U_{2}\left( \mathbb{F}_{q^{2}}\right) \right] $-submodule of $V_{k+(q-1)}$
and hence we get the well defined $U_{2}\left( \mathbb{F}_{q^{2}}\right) $%
-map of $\mathbb{F}_{q^{2}}$-spaces:%
\begin{equation*}
\overline{\theta }_{q}:V_{k-2}\rightarrow \dfrac{V_{k+(q-1)}}{D\left(
V_{k}\right) }\text{.}
\end{equation*}

\noindent Notice that this map is non-zero, then monic, i.e. $\func{Im}%
\theta _{q}\cap D\left( V_{k}\right) =0$, because our restriction on the
range of $k$ implies that $V_{k-2}^{^{\prime }}$ is an irreducible
representation of $G$.

By dimension considerations, we have the $\mathbb{F}_{q^{2}}$-vector space
decomposition $V_{k+(q-1)}=D\left( V_{k}\right) \oplus \theta _{q}\left(
V_{k-2}\right) \oplus W$, where%
\begin{equation}
W:=\dbigoplus_{r=k}^{q-1}\mathbb{F}_{q^{2}}X^{r}Y^{k+q-1-r}.  \label{WWW}
\end{equation}

Notice that the action of $\mathbb{F}_{q^{2}}^{\times }$ upon $W$ decomposes
the space into the direct sum of eigenspaces of dimension one: if $k\leq
r\leq q-1$, $\mathbb{F}_{q^{2}}X^{r}Y^{k+q-1-r}$ is the subspace of $W$ on
which $\mathbb{F}_{q^{2}}^{\times }$ acts by the character $t\mapsto
t^{r(q+1)+(q-1)-qk}$.

\begin{proposition}
\label{uffaaa}There exist $U_{2}\left( \mathbb{F}_{q^{2}}\right) $%
-epimorphisms of $\mathbb{F}_{q^{2}}$-spaces $\dfrac{V_{k+(q-1)}}{D\left(
V_{k}\right) }\rightarrow e^{k}\otimes V_{q-1-k}$ whose kernel is equal to $%
\func{Im}\overline{\theta }_{q}$ if and only if $q=p$ is prime; in this case
they are given by the maps $\left\{ \omega _{s}\right\} _{s\in \mathbb{F}%
_{q^{2}}^{\times }}$, where for any $s\in \mathbb{F}_{q^{2}}^{\times }$ the
homomorphism $\omega _{s}$ is defined by $\omega _{s}\left( \func{Im}%
\overline{\theta }_{q}\right) =0$ and:%
\begin{equation*}
X^{r}Y^{k+q-1-r}\func{mod}D\left( V_{k}\right) \mapsto s\cdot \dfrac{\tbinom{%
q-1-k}{r-k}}{\tbinom{q-1}{r}}\cdot X^{r-k}Y^{q-1-r},\text{ \ \ \ \ }\left(
k\leq r\leq q-1\right) .
\end{equation*}
\end{proposition}

\begin{proof}
We need to construct the $U_{2}\left( \mathbb{F}_{q^{2}}\right) $%
-equivariant onto maps $\omega :V_{k+(q-1)}\rightarrow e^{k}\otimes
V_{q-1-k} $ whose kernel is $D\left( V_{k}\right) \oplus \theta _{q}\left(
V_{k-2}\right) $, therefore we need to define $\omega $ only on $W$. Since $%
t\in \mathbb{F}_{q^{2}}^{\times }$ acts on the subspace $\mathbb{F}%
_{q^{2}}X^{s}Y^{q-1-k-s}\subset e^{k}\otimes V_{q-1-k}$ as multiplication by 
$t^{s(q+1)+(q-1)+kq}t^{k\left( 1-q\right) }$ ($0\leq s\leq q-1-k$),\ and
since $\omega $ has to respect the decomposition of the spaces in $\mathbb{F}%
_{q^{2}}^{\times }$-eigenspaces, we have $\omega :X^{r}Y^{k+q-1-r}\mapsto
\alpha \left( r\right) X^{r-k}Y^{q-1-r}$ for some $\alpha \left( r\right)
\in \mathbb{F}_{q^{2}}^{\times }$ ($k\leq r\leq q-1$). For any $k\leq r\leq
q-1$ we find the following conditions to be satisfied:%
\begin{eqnarray*}
\alpha \left( r\right) &=&\left( -1\right) ^{k}\alpha \left( q+k-1-r\right)
\\
\alpha \left( r\right) &=&\alpha \left( q-1\right) \cdot \dfrac{\tbinom{q-1-k%
}{r-k}}{\tbinom{q-1}{r}},
\end{eqnarray*}

where the first equation is equivalent to the equivariance of $\omega $ with
respect to the matrix $\left( 
\begin{array}{cc}
0 & 1 \\ 
-1 & 0%
\end{array}%
\right) $, and the second is equivalent to the equivariance with respect to $%
\left( 
\begin{array}{cc}
1 & 0 \\ 
u & 1%
\end{array}%
\right) $ for any $u\in \mathbb{F}$. Notice that $\tbinom{q-1-k}{r-k}$ is
not divisible by $p$ for every $k\leq r\leq q-1$ if and only if $q=p$, in
fact if $q>p$ we have $\tbinom{q-1-k}{p-k}\equiv 0(\func{mod}p).$

If $q=p$ we can choose $\alpha \left( q-1\right) \neq 0$ arbitrarily, so
that we get the maps in the statement; all of them are acceptable since:%
\begin{equation*}
U_{2}\left( \mathbb{F}_{q^{2}}\right) =\mathbb{F}_{q^{2}}^{\times }\cdot
\left\langle \left( 
\begin{array}{cc}
0 & 1 \\ 
-1 & 0%
\end{array}%
\right) ,\left( 
\begin{array}{cc}
1 & 0 \\ 
u & 1%
\end{array}%
\right) :u\in \mathbb{F}\right\rangle .
\end{equation*}
\end{proof}

Note that the above proposition gives, for $q>p$, equivariant non-surjective
maps $\dfrac{V_{k+(q-1)}}{D\left( V_{k}\right) }\rightarrow e^{k}\otimes
V_{q-1-k}$ whose kernel properly contains $\func{Im}\overline{\theta }_{q}$.

\begin{corollary}
\label{ES2}Let $q=p$ be prime. For any $s\in \mathbb{F}_{q^{2}}^{\times }$
we have an exact sequence of $\mathbb{F}_{q^{2}}\left[ U_{2}\left( \mathbb{F}%
_{q^{2}}\right) \right] $-modules:%
\begin{equation*}
0\rightarrow V_{k-2}\overset{\overline{\theta }_{q}}{\rightarrow }\dfrac{%
V_{k+(q-1)}}{D\left( V_{k}\right) }\overset{\omega _{s}}{\rightarrow }%
e^{k}\otimes V_{q-1-k}\rightarrow 0,
\end{equation*}

\noindent coming, if $s\in \mathbb{F}^{\times }$, from the exact sequence of 
$\mathbb{F}\left[ SL_{2}\left( \mathbb{F}\right) \right] $-modules:%
\begin{equation*}
0\rightarrow V_{k-2}^{\prime }\overset{\overline{\theta }_{q}}{\rightarrow }%
\dfrac{V_{k+(q-1)}^{\prime }}{D\left( V_{k}^{\prime }\right) }\overset{%
\omega _{s}}{\rightarrow }V_{q-1-k}^{\prime }\rightarrow 0.
\end{equation*}
\end{corollary}

\begin{remark}
If $q$ is any power of the prime $p$ and $2\leq k\leq p-1$, one sees
(formula \ref{ora et labora}) that there is an exact sequence of $\mathbb{F}%
\left[ GL_{2}\left( \mathbb{F}\right) \right] $-modules:%
\begin{equation*}
0\rightarrow e\otimes V_{k-2}^{\prime }\overset{\overline{\theta }_{q}}{%
\rightarrow }\dfrac{V_{k+(q-1)}^{\prime }}{D\left( V_{k}^{\prime }\right) }%
\rightarrow \limfunc{coker}\overline{\theta }_{q}\rightarrow 0.
\end{equation*}

\noindent The same reasoning of the proof of Proposition \ref{uffaaa} gives
that there exist $GL_{2}\left( \mathbb{F}\right) $-epimorphisms of $\mathbb{F%
}$-spaces $\frac{V_{k+(q-1)}^{\prime }}{D\left( V_{k}^{\prime }\right) }%
\rightarrow e^{k}\otimes V_{q-1-k}^{\prime }$ whose kernel is equal to $%
\func{Im}\overline{\theta }_{q}$ if and only if $q=p$ is prime (one just
needs to consider the equivariance of the desired maps with respect to the
characters of the split torus of $GL_{2}\left( \mathbb{F}\right) $, and
proceed as in the proof of the proposition).

\noindent Furthermore, a computation of Brauer characters gives the
following identity in $K_{0}\left( G,\mathbb{F}\right) $: 
\begin{equation*}
\limfunc{coker}\overline{\theta }_{q}=\left( V_{q-1-k}^{\prime }\right)
^{\ast }=e^{k}\cdot V_{q-1-k}^{^{\prime }}.
\end{equation*}

\noindent (notice that $(V_{k}^{\prime })^{\ast }=e^{-k}\cdot V_{k}^{\prime
} $ in $K_{0}\left( G,\mathbb{F}\right) $ for any $k\in 
\mathbb{Z}
$).\ \noindent In the case $q=p$, $e^{k}\otimes V_{q-1-k}^{^{\prime }}$ is
an irreducible $GL_{2}\left( \mathbb{F}\right) $-module, so that we find
(again) that $\limfunc{coker}\overline{\theta }_{q}\simeq e^{k}\otimes
V_{q-1-k}^{^{\prime }}$; if $q>p$ we cannot guarantee the irreducibility of $%
e^{k}\otimes V_{q-1-k}^{^{\prime }}$ and we can only deduce the existence of
isomorphisms:%
\begin{equation*}
(\limfunc{coker}\overline{\theta }_{q})^{ss}\simeq ((V_{q-1-k}^{\prime
})^{\ast })^{ss}\simeq (e^{k}\otimes V_{q-1-k}^{^{\prime }})^{ss}\text{.}
\end{equation*}
\end{remark}

\bigskip

We conclude this paragraph with the following result:

\begin{proposition}
\label{non-split}For any $2\leq k\leq p-1$,\ the exact sequence of $\mathbb{F%
}\left[ SL_{2}\left( \mathbb{F}\right) \right] $-modules $0\rightarrow
V_{k-2}^{\prime }\rightarrow \dfrac{V_{k+(q-1)}^{\prime }}{D\left(
V_{k}^{\prime }\right) }\rightarrow \limfunc{coker}\overline{\theta }%
_{q}\rightarrow 0$ is non-split.
\end{proposition}

\begin{proof}
If the sequence of $SL_{2}\left( \mathbb{F}\right) $-modules $0\rightarrow
V_{k-2}^{\prime }\rightarrow \frac{V_{k+(q-1)}^{\prime }}{D\left(
V_{k}^{\prime }\right) }\rightarrow \limfunc{coker}\overline{\theta }%
_{q}\rightarrow 0$ splits, then by tensoring with $\mathbb{F}_{q^{2}}$ we
would get a split sequence of $U_{2}\left( \mathbb{F}_{q^{2}}\right) $%
-modules, since $SL_{2}\left( \mathbb{F}\right) \lhd U_{2}\left( \mathbb{F}%
_{q^{2}}\right) $ is of index prime to $p$.

\noindent Let us assume that $0\rightarrow V_{k-2}\overset{\overline{\theta }%
_{q}}{\rightarrow }\frac{V_{k+(q-1)}}{D\left( V_{k}\right) }\rightarrow 
\limfunc{coker}\overline{\theta }_{q}\rightarrow 0$ is a split exact
sequence of $\mathbb{F}_{q^{2}}\left[ U_{2}\left( \mathbb{F}_{q^{2}}\right) %
\right] $-modules. A basis of $\frac{V_{k+(q-1)}}{D\left( V_{k}\right) }$ is
given by the set:%
\begin{equation*}
\mathcal{S}=\left\{ X^{r}Y^{k+q-1-r}\func{mod}D(V_{k}),X^{s}Y^{k-2-s}\theta
_{q}\func{mod}D(V_{k})\right\} _{k\leq r\leq q-1,0\leq s\leq k-2},
\end{equation*}

and the group $\mathbb{F}_{q^{2}}^{\times }$ acts on the $1$-dimensional
space generated by each of the vectors in $\mathcal{S}$ as multiplication by
a character, and this action is via different characters on different spaces
(in fact, $\mathbb{F}_{q^{2}}^{\times }$ acts on $X^{s}Y^{k-2-s}\theta _{q}%
\func{mod}D(V_{k})$ via the character $t\mapsto t^{s(q+1)-kq+2q}$ for any $%
0\leq s\leq k-2$; for the action on the other basis vectors see the
computation before Proposition \ref{uffaaa}). This implies that if $M$ is a
complement of $\func{Im}\overline{\theta }_{q}$ inside $\frac{V_{k+(q-1)}}{%
D\left( V_{k}\right) }$, there is basis of $M$ contained in $\mathcal{S}$;
since $M\cap \func{Im}\overline{\theta }_{q}=0$ we have $M=(W\oplus D\left(
V_{k}\right) )/D\left( V_{k}\right) $ (see \ref{WWW} for the definition of $%
W $). This is not possible, since this space is not closed under the action
of $U_{2}\left( \mathbb{F}_{q^{2}}\right) $.\ In fact we have:%
\begin{equation*}
\begin{pmatrix}
1 & 1 \\ 
1 & 0%
\end{pmatrix}%
\cdot X^{k}Y^{q-1}\func{mod}D(V_{k})=\sum\nolimits_{i=0}^{k}\binom{k}{i}%
X^{i}Y^{k+q-1-i}\func{mod}D(V_{k}),
\end{equation*}

and $\sum\nolimits_{i=0}^{k}\binom{k}{i}X^{i}Y^{k+q-1-i}=X^{k}Y^{q-1}+\tsum%
\nolimits_{i=0}^{k-1}\binom{k}{i}X^{i}Y^{k+q-1-i}$, where the first term is
in $W$ and the second is not in $W\oplus D\left( V_{k}\right) $.
\end{proof}

\begin{corollary}
If $2\leq k\leq p-1$ and $q=p$ is prime, the one dimensional $\mathbb{F}$%
-vector space $\limfunc{Ext}_{SL_{2}\left( \mathbb{F}\right) }^{1}\left(
V_{k-2}^{\prime },V_{q-1-k}^{\prime }\right) \simeq \limfunc{Ext}%
_{SL_{2}\left( \mathbb{F}\right) }^{1}\left( V_{q-1-k}^{\prime
},V_{k-2}^{\prime }\right) $ is generated by the class of $\dfrac{V_{k+(q-1)}%
}{D\left( V_{k}\right) }$.
\end{corollary}

\begin{proof}
This follows from the previous result and from the Propositions in \S \S %
12.1,12.2 of \cite{Hum}.
\end{proof}

\subsection{Connecting the two exact sequences}

We would like to be able to relate the exact sequences coming from Corollary %
\ref{ES1} (for $\mathbb{K=F}_{q^{2}}$) and Corollary \ref{ES2} (with the
additional assumption $q=p$), by constructing - for some integers $\alpha
,\beta ,\gamma $ - a commutative diagram of the following form:%
\begin{equation*}
\begin{array}{ccccccccc}
0 & \rightarrow & H^{0}(C,\Omega _{C}^{1})_{-k} & \rightarrow & 
H_{dR}^{1}(C)_{-k} & \rightarrow & H^{1}(C,\mathcal{O}_{C})_{-k} & 
\rightarrow & 0 \\ 
&  & \downarrow \simeq &  & \downarrow \simeq &  & \downarrow \simeq &  & 
\\ 
0 & \rightarrow & e^{\alpha }\otimes V_{k-2} & \rightarrow & e^{\beta
}\otimes \dfrac{V_{k+q-1}}{D\left( V_{k}\right) } & \rightarrow & e^{\gamma
}\otimes V_{q-1-k} & \rightarrow & 0.%
\end{array}%
\end{equation*}

\noindent In the sequel, we will construct the vertical isomorphisms; we
will also treat the case $q>p$.

\subsubsection{The left vertical maps}

We assume from now on $\mathbb{K=F}_{q^{2}}$ and $C$ will denote the smooth
projective curve $C_{/\mathbb{F}_{q^{2}}}$; as usual we fix the integer $k$
such that $2\leq k\leq p-1$ ($p>2$). We will always allow $q$ to be any
positive power of $p$, unless otherwise stated.

There is a natural isomorphism of $U_{2}\left( \mathbb{F}_{q^{2}}\right) $%
-modules $e^{-1}\otimes H^{0}\left( \mathbb{P}^{2}(\mathbb{F}_{q^{2}}),%
\mathcal{O}\left( q-2\right) \right) \simeq e^{-1}\otimes \limfunc{Sym}%
^{q-2}(\mathbb{F}_{q^{2}}^{3})$; furthermore by \cite{HaaJan}, Prop. 2.1.
the map:

\begin{eqnarray*}
e^{-1}\otimes H^{0}\left( \mathbb{P}^{2}(\mathbb{F}_{q^{2}}),\mathcal{O}%
\left( q-2\right) \right) &\rightarrow &H^{0}\left( C,\Omega _{C}^{1}\right)
\\
X^{b}Y^{a}Z^{q-2-a-b} &\mapsto &-t_{1}^{a}t_{2}^{q-2-a-b}dt_{2},
\end{eqnarray*}%
where $a,b\geq 0,a+b\leq q-2$ and $t_{1}=Y/X,t_{2}=Z/X$, is a \textit{natural%
} $U_{2}\left( \mathbb{F}_{q^{2}}\right) $-isomorphism. Here the action of $%
U_{2}\left( \mathbb{F}_{q^{2}}\right) $ on the codomain is given by
embedding $U_{2}\left( \mathbb{F}_{q^{2}}\right) $ in $GL_{3}\left( \mathbb{F%
}_{q^{2}}\right) $ via $g\mapsto \left( 
\begin{array}{cc}
(g^{-1})^{t} & 0 \\ 
0 & 1%
\end{array}%
\right) $ and letting the latter group act on $\mathbb{F}_{q^{2}}\left[ X,Y,Z%
\right] $ in the standard way (see formula \ref{action} in Section 2). For
example we have, for $a,b\geq 0$ and $a+b\leq q-2$:

\begin{eqnarray}
\left( 
\begin{array}{cc}
0 & 1 \\ 
-1 & 0%
\end{array}%
\right) \cdot t_{1}^{a}t_{2}^{b}dt_{2} &=&\left( 
\begin{array}{ccc}
0 & 1 & 0 \\ 
-1 & 0 & 0 \\ 
0 & 0 & 1%
\end{array}%
\right) \cdot \left( \frac{Y}{X}\right) ^{a}\left( \frac{Z}{X}\right)
^{b}d\left( \frac{Z}{X}\right) =  \notag \\
&=&\left( \frac{X}{-Y}\right) ^{a}\left( \frac{Z}{-Y}\right) ^{b}d\left( 
\frac{Z}{-Y}\right)  \notag \\
&=&\left( -t_{1}^{-1}\right) ^{a}\left( -t_{1}^{-1}t_{2}\right) ^{b}d\left(
-t_{1}^{-1}t_{2}\right) ;
\end{eqnarray}%
then we use the identity $t_{2}^{q+1}=t_{1}^{q}-t_{1}$ to get $%
dt_{1}=-t_{2}^{q}dt_{2},d(t_{1}^{-1}t_{2})=t_{1}^{q-2}dt_{2}$ and conclude
that:

\begin{equation*}
\left( 
\begin{array}{cc}
0 & 1 \\ 
-1 & 0%
\end{array}%
\right) t_{1}^{a}t_{2}^{b}dt_{2}=\left( -1\right)
^{a+b+1}t_{1}^{q-2-a-b}t_{2}^{b}dt_{2}.
\end{equation*}

\noindent We recall for the sequel the action of the group $U_{2}\left( 
\mathbb{F}_{q^{2}}\right) $ on every element (for a proof, cf. \cite{HaaJan}%
): for any $a,b\geq 0$ and $a+b\leq q-2$:

\begin{equation}
\left( 
\begin{array}{cc}
0 & 1 \\ 
-1 & 0%
\end{array}%
\right) t_{1}^{a}t_{2}^{b}dt_{2}=\left( -1\right)
^{a+b+1}t_{1}^{q-2-a-b}t_{2}^{b}dt_{2}  \label{1g}
\end{equation}

\begin{equation}
\left( 
\begin{array}{cc}
1 & 0 \\ 
u & 1%
\end{array}%
\right) t_{1}^{a}t_{2}^{b}dt_{2}=\tsum\nolimits_{i=0}^{a}\tbinom{a}{i}\left(
-u\right) ^{a-i}t_{1}^{i}t_{2}^{b}dt_{2},\text{ }\ (u\in \mathbb{F)}
\label{2g}
\end{equation}

\begin{equation}
\left( 
\begin{array}{cc}
t & 0 \\ 
0 & t^{-q}%
\end{array}%
\right) t_{1}^{a}t_{2}^{b}dt_{2}=t^{a(q+1)+b+1}t_{1}^{a}t_{2}^{b}dt_{2},%
\text{ \ \ \ \ \ \ \ \ \ }\ \ (t\in \mathbb{F}_{q^{2}}^{\times }).
\label{3g}
\end{equation}

\noindent Therefore we can make the identification $H^{0}\left( C,\Omega
_{C}^{1}\right) _{-k}=\bigoplus\nolimits_{a=0}^{k-2}\mathbb{F}%
_{q^{2}}t_{1}^{a}t_{2}^{q-k}dt_{2}$.

\begin{proposition}
\label{prop}The isomorphisms of $\mathbb{F}_{q^{2}}[U_{2}\left( \mathbb{F}%
_{q^{2}}\right) ]$-modules $H^{0}\left( C,\Omega _{C}^{1}\right)
_{-k}\rightarrow e^{1-k}\otimes V_{k-2}$ are the maps $\left\{ \varphi
_{s}\right\} _{s\in \mathbb{F}_{q^{2}}^{\times }}$, where for any $s\in 
\mathbb{F}_{q^{2}}^{\times }$, $\varphi _{s}$ is defined by:%
\begin{equation*}
t_{1}^{a}t_{2}^{q-k}dt_{2}\mapsto \left( -1\right) ^{a}\cdot s\cdot
X^{a}Y^{k-2-a},\text{ \ \ \ \ }\left( 0\leq a\leq k-2\right) .
\end{equation*}
\end{proposition}

\begin{proof}
The isomorphisms $\varphi $ we are looking for preserve the $\mathbb{F}%
_{q^{2}}^{\times }$-eigenspace decomposition of the modules: for any fixed $%
0\leq i\leq k-2$, $\mathbb{F}_{q^{2}}X^{i}Y^{k-2-i}\subset e^{1-k}\otimes
V_{k-2}$ is the eigenspace relative to the character of $\mathbb{F}%
_{q^{2}}^{\times }$ represented by the $\func{mod}\left( q^{2}-1\right) $
integer $(i+1)(q+1)-k$; $\mathbb{F}_{q^{2}}t_{1}^{a}t_{2}^{q-k}dt_{2}\subset
H^{0}\left( C,\Omega _{C}^{1}\right) _{-k}$ is the eigenspace for $%
a(q+1)+q-k+1$ ($0\leq a\leq k-2$). We deduce that $\varphi $ has to send $%
t_{1}^{a}t_{2}^{q-k}dt_{2}$ to $\beta \left( a\right) X^{a}Y^{k-2-a}$ for
any $0\leq a\leq k-2$; by imposing the condition of equivariance with
respect to the matrices $\left( 
\begin{array}{cc}
0 & 1 \\ 
-1 & 0%
\end{array}%
\right) ,\left( 
\begin{array}{cc}
1 & 0 \\ 
u & 1%
\end{array}%
\right) $ ($u\in \mathbb{F}$), we find respectively, for any $0\leq a\leq
k-2 $:%
\begin{eqnarray*}
\beta \left( a\right) &=&\left( -1\right) ^{k}\beta \left( k-2-a\right) \\
\beta \left( a\right) &=&\left( -1\right) ^{k+a}\beta \left( k-2\right) .
\end{eqnarray*}

\noindent Any fixed value for $\beta \left( k-2\right) \in \mathbb{F}%
_{q^{2}}^{\times }$ makes these equations satisfied.
\end{proof}

\subsubsection{The central vertical maps}

We recall a computation from \cite{HaaJan}, \S 4: let $U_{0}=\limfunc{Spec}%
\mathbb{F}_{q^{2}}[Y/X,Z/X]$, $U_{1}=\limfunc{Spec}\mathbb{F}%
_{q^{2}}[X/Y,Z/Y]$ be open in $\mathbb{P}^{2}\left( \mathbb{F}%
_{q^{2}}\right) $ and let $\mathfrak{U}=\{U_{0}\cap C,U_{1}\cap C\}.$ Using 
\v{C}ech cohomology with respect to open affine covering $\mathfrak{U}$ of $%
C $, we recover $H_{dR}^{1}(C)$ as the quotient of the space:%
\begin{equation*}
\left\{ \left( \omega _{0},\omega _{1},f_{01}\right) :\omega _{i}\in
H^{0}(U_{i}\cap C,\Omega _{C}^{1}),f_{01}\in H^{0}\left( U_{i}\cap C,%
\mathcal{O}_{C}\right) ,df_{01}=\omega _{0}-\omega _{1}\right\} ,
\end{equation*}

\noindent by the subspace $\left\{ \left( df_{0},df_{1},f_{0}-f_{1}\right)
:f_{i}\in H^{0}(U_{i}\cap C,\mathcal{O}_{C})\right\} $.

\noindent Furthermore, the non-trivial maps in the exact sequence in
Proposition \ref{here} are induced by $\omega \mapsto \left[ \left( \omega
_{|U_{0}\cap C},\omega _{|U_{1}\cap C},0\right) \right] $ and by $\left[
\left( \omega _{0},\omega _{1},f_{01}\right) \right] \mapsto f_{01}$.

Once we fix the integer $2\leq k\leq p-1$ and we pass to the $-k$
eigenspaces of the above cohomology groups, we see that a basis\footnote{%
The basis constructed in \cite{HaaJan} works in general for any $1\leq k\leq
q$.} for $H_{dR}^{1}(C)_{-k}$ is built up by two sets: $\mathcal{A}=\left\{
e_{a,q-k}\right\} _{a=0}^{k-2}$ and $\mathcal{B}=\left\{ e_{a,q+1-k}\right\}
_{-a=1}^{q-k}$, where:%
\begin{equation*}
e_{a,q-k}:=[(t_{1}^{a}t_{2}^{q-k}dt_{2},t_{1}^{a}t_{2}^{q-k}dt_{2},0)]
\end{equation*}

\begin{equation*}
e_{a,q+1-k}:=[(-at_{1}^{q-1+a}t_{2}^{q-k}dt_{2},-(a+1-k)t_{1}^{a}t_{2}^{q-k}dt_{2},t_{1}^{a}t_{2}^{q+1-k})].
\end{equation*}

\noindent This follows from the previous observations and from Proposition
2.3. in \cite{HaaJan}, that gives explicitly a basis for $H^{1}(C,\mathcal{O}%
_{C})$.

Notice that the action of $U_{2}\left( \mathbb{F}_{q^{2}}\right) $ on $%
\mathcal{A}$ is clear (see formulae \ref{1g}, \ref{2g} and \ref{3g}), while
for $\mathcal{B}$ we refer to \cite{HaaJan} for the following computation:
for any integer $a$ such that $1\leq -a\leq q-k$ and any $u\in \mathbb{F}$, $%
t\in \mathbb{F}_{q^{2}}^{\times }$, one has:

\begin{equation}
\left( 
\begin{array}{cc}
0 & 1 \\ 
-1 & 0%
\end{array}%
\right) e_{a,q+1-k}=\left( -1\right) ^{a+k+1}e_{-(a+q+1-k),q+1-k},  \label{4}
\end{equation}

\begin{center}
\begin{equation}
\begin{tabular}{l}
$\left( 
\begin{array}{cc}
1 & 0 \\ 
u & 1%
\end{array}%
\right) e_{a,q+1-k}=\tsum\nolimits_{i=0}^{a+q-k}\tbinom{q+a}{i}\left(
-u\right) ^{i}e_{a-i,q+1-k}+$ \\ 
$\ \ \ \ \ \ \ \ \ \ \ \ \ \ \ \ \ \ \ \ \ \ \ \ \ \ \ \ \ \
+\tsum\nolimits_{i=a+q+1-k}^{a+q-1}\tbinom{q+a-1}{i}\left( -u\right)
^{i}ae_{q+a-1-i,q-k},$%
\end{tabular}
\label{5}
\end{equation}
\end{center}

\begin{equation}
\left( 
\begin{array}{cc}
t & 0 \\ 
0 & t^{-q}%
\end{array}%
\right) e_{a,q+1-k}=t^{a(q+1)+q-k+1}e_{a,q+1-k}.  \label{6}
\end{equation}

Observe that the canonical projection $\pi :H_{dR}^{1}(C)_{-k}\rightarrow
H^{1}(C,\mathcal{O})_{-k}$ coming from Corollary \ref{ES1} is given by
sending the elements of $\mathcal{A}$ to zero and by $e_{a,q+1-k}\mapsto
t_{1}^{a}t_{2}^{q+1-k}$, for any $1\leq -a\leq q-k$.

We can now prove:

\begin{proposition}
\label{long}The $\mathbb{F}_{q^{2}}\left[ U_{2}\left( \mathbb{F}%
_{q^{2}}\right) \right] $ modules $H_{dR}^{1}\left( C\right) _{-k}$ and $%
e^{1-k}\otimes \dfrac{V_{k+(q-1)}}{D\left( V_{k}\right) }$ are isomorphic. A
family of $U_{2}\left( \mathbb{F}_{q^{2}}\right) $-equivariant isomorphisms
is given by the maps $\left\{ f_{s}\right\} _{s\in \mathbb{F}%
_{q^{2}}^{\times }}$, where, for any $s\in \mathbb{F}_{q^{2}}^{\times }$, $%
f_{s}$ is defined by:%
\begin{equation*}
f_{s}\left( e_{a,q-k}\right) =\left( -1\right) ^{a}\cdot s\cdot
X^{a}Y^{k-2-a}\theta _{q}\func{mod}D\left( V_{k}\right) ,\text{ \ \ \ \ \ \
for\ }0\leq a\leq k-2,
\end{equation*}%
\begin{equation*}
f_{s}\left( e_{a,q+1-k}\right) =\left( -1\right) ^{a}\cdot k\cdot s\cdot
X^{a+q}Y^{k-1-a}\func{mod}D\left( V_{k}\right) ,\text{ \ \ for }1\leq -a\leq
q-k.
\end{equation*}
\end{proposition}

\begin{proof}
Let $f$ be an isomorphism as above; the action of $\mathbb{F}%
_{q^{2}}^{\times }$ decomposes both $H_{dR}^{1}\left( C\right) _{-k}$ and $%
e^{1-k}\otimes \frac{V_{k+(q-1)}}{D\left( V_{k}\right) }$ in the direct sum
of \textit{one dimensional} eigenspaces, therefore an easy computation shows
that:%
\begin{eqnarray*}
f\left( e_{a,q-k}\right) &=&\gamma \left( a\right) \cdot
X^{a}Y^{k-2-a}\theta _{q}\func{mod}D\left( V_{k}\right) \text{ \ for }0\leq
a\leq k-2, \\
f\left( e_{a,q+1-k}\right) &=&\delta \left( a\right) \cdot X^{a+q}Y^{k-1-a}%
\func{mod}D\left( V_{k}\right) \text{ \ \ for }1\leq -a\leq q-k.
\end{eqnarray*}

\noindent where $\gamma \left( a\right) ,\delta \left( a\right) \in \mathbb{F%
}_{q^{2}}^{\times }$. Since the span of the $e_{a,q-k}$'s for $0\leq a\leq
k-2$ is just $H^{0}(C,\Omega _{C}^{1})_{-k}$ and it has to be mapped by $f$\
isomorphically onto $\overline{\theta }_{q}\left( e^{1-k}\otimes
V_{k-2}\right) \simeq e^{1-k}\otimes V_{k-2}$, Proposition \ref{prop} gives $%
\gamma \left( a\right) =\left( -1\right) ^{a}s$ for some $s\in \mathbb{F}%
_{q^{2}}^{\times }$ that we consider now fixed.

Let $1\leq -a\leq q-k$;\ the relation: 
\begin{equation*}
f\left( \left( 
\begin{array}{cc}
0 & 1 \\ 
-1 & 0%
\end{array}%
\right) e_{a,q+1-k}\right) =\left( 
\begin{array}{cc}
0 & 1 \\ 
-1 & 0%
\end{array}%
\right) \cdot \delta \left( a\right) X^{a+q}Y^{k-1-a}\func{mod}D\left(
V_{k}\right)
\end{equation*}
gives $\left[ (-1)^{a+q}\delta \left( a\right) -\left( -1\right)
^{a+k+1}\delta \left( k-a-q-1\right) \right] X^{k-1-a}Y^{a+q}\in D\left(
V_{k}\right) $; by a degree consideration this can only happen when the
coefficient of the monomial is zero, that is when:%
\begin{equation*}
\delta \left( a\right) =\left( -1\right) ^{k}\delta \left( k-a-q-1\right) .
\end{equation*}

\noindent For any fixed $u\in \mathbb{F}$ we also require: 
\begin{equation*}
f\left( \left( 
\begin{array}{cc}
1 & 0 \\ 
u & 1%
\end{array}%
\right) e_{a,q+1-k}\right) =\left( 
\begin{array}{cc}
1 & 0 \\ 
u & 1%
\end{array}%
\right) \cdot \delta \left( a\right) X^{a+q}Y^{k-1-a}\func{mod}D\left(
V_{k}\right) .
\end{equation*}%
By making the action of the matrix explicit, the condition we find after
some simplifications is:

\medskip

$\ \tsum\nolimits_{i=0}^{a+q-k}\tbinom{q+a}{i}\left( -u\right) ^{i}\delta
\left( a-i\right) X^{a-i+q}Y^{k-1-a+i}+$

\medskip

$+\tsum\nolimits_{i=a+q-k+1}^{a+q-1}\tbinom{q+a}{i}\left( a-i\right) \left(
-u\right) ^{i}\gamma \left( q+a-1-i\right) X^{a-i+2q-1}Y^{k-q-a+i}+$

\medskip

$-\tsum\nolimits_{i=a+q-k+1}^{a+q-1}\tbinom{q+a}{i}\left( a-i\right) \left(
-u\right) ^{i}\gamma \left( q+a-1-i\right) X^{a-i+q}Y^{k-1-a+i}\equiv $

\medskip

$\equiv \tsum\nolimits_{i=0}^{a+q-k}\tbinom{q+a}{i}u^{i}\delta \left(
a\right) X^{a-i+q}Y^{k-1-a+i}+$

\medskip

$+\tsum\nolimits_{i=a+q-k+1}^{a+q-1}\tbinom{q+a}{i}u^{i}\delta \left(
a\right) X^{a-i+q}Y^{k-1-a+i}$ $\ \func{mod}D\left( V_{k}\right) .$

\medskip

\noindent This is equivalent to:

\medskip

$\tsum\nolimits_{i=0}^{a+q-k}\tbinom{q+a}{i}u^{i}\left[ \left( -1\right)
^{i}\delta \left( a-i\right) -\delta \left( a\right) \right]
X^{a-i+q}Y^{k-1-a+i}+$

\medskip

$+\tsum\nolimits_{i=a+q-k+1}^{a+q-1}\tbinom{q+a}{i}\left( a-i\right)
u^{i}\left( -1\right) ^{a}sX^{a-i+2q-1}Y^{k-q-a+i}+$

\medskip

$-\tsum\nolimits_{i=a+q-k+1}^{a+q-1}\tbinom{q+a}{i}u^{i}\left[ \left(
a-i\right) \left( -1\right) ^{a}s+\delta \left( a\right) \right]
X^{a-i+q}Y^{k-1-a+i}\equiv 0\ \func{mod}D\left( V_{k}\right) .$

\medskip

Notice that for $0\leq i\leq a+q-k$, the monomial $X^{a-i+q}Y^{k-1-a+i}$
does not belong to $D\left( V_{k}\right) $.\ The above congruence holds 
\textit{if} $\delta \left( a\right) =\left( -1\right) ^{i}\delta \left(
a-i\right) $ (for any $i=0,...,a+q-k$ and any $-a=1,...,q-k$) and the
polynomial $\tsum\nolimits_{i=a+q-k+1}^{a+q-1}\tbinom{q+a}{i}u^{i}\cdot
g_{a,i,s}\left( X,Y\right) $ is in $D(V_{k})$, where:$\ $%
\begin{eqnarray*}
g_{a,i,s}\left( X,Y\right) &=&\left( a-i\right) \left( -1\right)
^{a}sX^{\left( a-i+q\right) +(q-1)}Y^{k-q-a+i}+ \\
&&-\left[ \left( a-i\right) \left( -1\right) ^{a}s+\delta \left( a\right) %
\right] X^{a-i+q}Y^{\left( k-q-a+i\right) +(q-1)}.
\end{eqnarray*}

In our range, we have $1\leq a+q-i\leq k-1$; since a basis for $D\left(
V_{k}\right) $ is give by: 
\begin{equation*}
\left\{ jX^{j+q-1}Y^{k-j}+(k-j)X^{j}Y^{k-j+q-1}\right\} _{j=0}^{k},
\end{equation*}%
we deduce that the above polynomial is in $D\left( V_{k}\right) $ if \ for
any $a$ such that $1\leq -a\leq q-k$ we can find an element $\sigma _{a}\in 
\mathbb{F}_{q^{2}}^{\times }$ such that $\sigma _{a}\cdot \left( a-i\right)
=\left( -1\right) ^{a}\left( a-i\right) s$ and $\sigma _{a}\cdot \left(
k+i-a\right) =-\left( -1\right) ^{a}\left( a-i\right) s-\delta \left(
a\right) $. This is equivalent to $\delta \left( a\right) =ks\left(
-1\right) ^{a+1}$ for every $1\leq -a\leq q-k$. This expression for $\delta
\left( \cdot \right) $ is compatible with all the previous conditions we
imposed on it.
\end{proof}

\subsubsection{The right vertical maps}

We now complete the construction of the main diagram. Recall that $2\leq
k\leq p-1$ and denote by $\pi ^{\prime }$ (resp. $\pi ^{\prime \prime }$)
the canonical projection $H_{dR}^{1}\left( C\right) _{-k}\rightarrow 
\limfunc{coker}\iota $ (resp. $e^{1-k}\otimes \frac{V_{k+(q-1)}}{D\left(
V_{k}\right) }\rightarrow \limfunc{coker}\overline{\theta }_{q}$), where $%
\iota :H^{0}\left( C,\Omega _{C}^{1}\right) _{-k}\hookrightarrow
H_{dR}^{1}\left( C\right) _{-k}$. We have:

\begin{corollary}
For any $s\in \mathbb{F}_{q^{2}}^{\times }$ there is a unique $U_{2}\left( 
\mathbb{F}_{q^{2}}\right) $-equivariant isomorphism of $\mathbb{F}_{q^{2}}$%
-spaces $\psi _{s}:\limfunc{coker}\iota \rightarrow \limfunc{coker}\overline{%
\theta }_{q}$ such that $\psi _{s}\circ \pi ^{\prime }=\pi ^{\prime \prime
}\circ f_{s}$. In particular we have $\limfunc{coker}\overline{\theta }%
_{q}\simeq H^{1}(C,\mathcal{O}_{C})_{-k}.$
\end{corollary}

\begin{proof}
By Proposition \ref{prop} and Proposition \ref{long} we have the commutative
diagram:%
\begin{equation*}
\begin{array}{ccccc}
0 & \rightarrow & H^{0}(C,\Omega _{C}^{1})_{-k} & \overset{\iota }{%
\rightarrow } & H_{dR}^{1}(C)_{-k} \\ 
&  & \downarrow \varphi _{s} &  & \downarrow f_{s} \\ 
0 & \rightarrow & e^{1-k}\otimes V_{k-2} & \overset{\overline{\theta }_{q}}{%
\rightarrow } & e^{1-k}\otimes \dfrac{V_{k+(q-1)}}{D\left( V_{k}\right) }%
\end{array}%
\end{equation*}

\noindent whose rows are exact. The statement just follows from this and the
identification $\limfunc{coker}\iota \simeq H^{1}(C,\mathcal{O}_{C})_{-k}.$
\end{proof}

\begin{remark}
If $q=p$, by Corollary \ref{ES2}, we have an isomorphism of $U_{2}\left( 
\mathbb{F}_{q^{2}}\right) $-modules $H^{1}(C,\mathcal{O}_{C})_{-k}\simeq
e\otimes V_{q-1-k}.$ By \cite{HaaJan}, Prop. 2.8. we also have, for any $q$,
the isomorphisms of $SL_{2}\left( \mathbb{F}\right) $-modules:%
\begin{equation*}
\limfunc{coker}\overline{\theta }_{q}\simeq H^{1}(C,\mathcal{O}%
_{C})_{-k}\simeq V_{q-1-k}^{\ast }.
\end{equation*}

Notice that if $q=p$ then $V_{q-1-k}$ is irreducible as representation of $%
SL_{2}\left( \mathbb{F}\right) $, and self-dual.
\end{remark}

\subsubsection{Back to the base field}

Now we go back to our original base field $\mathbb{F=F}_{q}$, $q=p^{n}$. At
this purpose, we consider the projective curve $C$ as defined over $\mathbb{F%
}$, so that its various cohomology groups are $\mathbb{F}$-spaces. Notice
that the elements of the basis we considered for the cohomology groups of $%
C_{/\mathbb{F}_{q^{2}}}$ are defined over the prime field of $\mathbb{F}$,
so that we will write $H_{dR}^{1}\left( C_{/\mathbb{F}}\right) _{-k}$ to
denote the $\mathbb{F}$-space with basis $\left\{ e_{a,q-k}\right\}
_{a=0}^{k-2}\cup \left\{ e_{a,q+1-k}\right\} _{-a=1}^{q-k}$ and similarly
for the other cohomology (sub)spaces we have been treating. In particular
the sequence of $\mathbb{F}\left[ SL_{2}\left( \mathbb{F}\right) \right] $%
-modules:\ 
\begin{equation*}
\begin{array}{ccccccccc}
0 & \rightarrow & H^{0}(C_{/\mathbb{F}},\Omega _{C_{/\mathbb{F}}}^{1})_{-k}
& \overset{\iota }{\rightarrow } & H_{dR}^{1}(C_{/\mathbb{F}})_{-k} & 
\overset{\pi }{\rightarrow } & H^{1}(C_{/\mathbb{F}},\mathcal{O}_{C_{/%
\mathbb{F}}})_{-k} & \rightarrow & 0%
\end{array}%
\end{equation*}

\noindent remains exact. (Notice that we do not change the names of the maps
involved with respect to the case of modules over $\mathbb{F}_{q^{2}}$ if
the behavior of the map is clear as above).

We can summarize some of the results we obtained as follows:

\begin{proposition}
\label{Main-SL2}Let $q>2$ be a power of a rational prime $p$, $\mathbb{F}$ a
field with $q$ elements, and let $k$ be an integer such that $2\leq k\leq
p-1 $. Let $SL_{2}\left( \mathbb{F}\right) $ act on $\mathbb{F}^{2}$ in the
standard way, and define the left $SL_{2}\left( \mathbb{F}\right) $%
-representation $V_{k}^{\prime }:=\limfunc{Sym}^{k}\mathbb{F}^{2}$. Define
furthermore $C$ to be the projective curve associated to the polynomial $%
XY^{q}-X^{q}Y-Z^{q+1}\in \mathbb{F}\left[ X,Y,Z\right] $.\ For any $s\in 
\mathbb{F}^{\times }$ we have the following commutative diagram of $\mathbb{F%
}\left[ SL_{2}\left( \mathbb{F}\right) \right] $-modules, in which the rows
are exact and the vertical arrows are isomorphisms:%
\begin{equation*}
\begin{array}{ccccccccc}
0 & \rightarrow & H^{0}(C_{/\mathbb{F}},\Omega _{C_{/\mathbb{F}}}^{1})_{-k}
& \overset{\iota }{\rightarrow } & H_{dR}^{1}(C_{/\mathbb{F}})_{-k} & 
\overset{\pi }{\rightarrow } & H^{1}(C_{/\mathbb{F}},\mathcal{O}_{C_{/%
\mathbb{F}}})_{-k} & \rightarrow & 0 \\ 
&  & \downarrow \varphi _{s} &  & \downarrow f_{s} &  & \downarrow \psi _{s}
&  &  \\ 
0 & \rightarrow & e^{1-k}\otimes V_{k-2}^{\prime } & \overset{\overline{%
\theta }_{q}}{\rightarrow } & e^{1-k}\otimes \dfrac{V_{k+(q-1)}^{\prime }}{%
D\left( V_{k}^{\prime }\right) } & \overset{\pi ^{\prime \prime }}{%
\rightarrow } & \limfunc{coker}\overline{\theta }_{q} & \rightarrow & 0.%
\end{array}%
\end{equation*}
\end{proposition}

\section{Integral models}

From now on we work over the field $\mathbb{F}$; we then go back to the
notation of Section 2 and write $V_{h}$ for the left $G=GL_{2}\left( \mathbb{%
F}\right) $-representation $\limfunc{Sym}^{h}\mathbb{F}^{2}$ ($h\geq 0$).

\subsection{Extending the action to $GL_{2}$}

Let $k$ be an integer such that $2\leq k\leq p-1$; the natural action of $G$
on $V_{k-2}$ and on $\frac{V_{k+(q-1)}}{D\left( V_{k}\right) }$ gives an
exact sequence of $\mathbb{F}\left[ G\right] $-modules:

\begin{equation}
0\rightarrow e\otimes V_{k-2}\overset{\overline{\theta }_{q}}{\rightarrow }%
\dfrac{V_{k+(q-1)}}{D\left( V_{k}\right) }\overset{\pi ^{\prime \prime }}{%
\rightarrow }\limfunc{coker}\overline{\theta }_{q}\rightarrow 0.
\label{ora et labora}
\end{equation}

\noindent For any $s\in \mathbb{F}^{\times }$, there is a unique extension
of the action of $SL_{2}\left( \mathbb{F}\right) $ on $H^{0}(C_{/\mathbb{F}%
},\Omega _{C_{/\mathbb{F}}}^{1})_{-k}$ and on $H_{dR}^{1}(C_{/\mathbb{F}%
})_{-k}$ to an action of $G$ making $\varphi _{s}$ and $f_{s}$ isomorphisms
of $G$-modules. By definition of $\varphi _{s}$ and $f_{s}$ we see that this
extension does not depend upon $s$ and is given by defining, for any $%
x_{1},x_{2}\in \mathbb{F}^{\times }:$

\begin{eqnarray*}
\left( 
\begin{array}{cc}
x_{1} & 0 \\ 
0 & x_{2}%
\end{array}%
\right) \cdot e_{a,q-k} &:&=x_{1}^{a+1}x_{2}^{k-a-1}e_{a,q-k}\text{, \ \ \ \
\ \ \ \ }(0\leq a\leq k-2) \\
\text{ }\left( 
\begin{array}{cc}
x_{1} & 0 \\ 
0 & x_{2}%
\end{array}%
\right) \cdot e_{a,q+1-k} &:&=x_{1}^{a+1}x_{2}^{k-a-1}e_{a,q+1-k}\text{, \ \
\ \ }(1\leq -a\leq q-k).
\end{eqnarray*}

Note that this definition of the action of the standard maximal torus of $G$
on $H_{dR}^{1}(C_{/\mathbb{F}})_{-k}$ is forced if we want it to be
compatible with any $U_{2}\left( \mathbb{F}_{q^{2}}\right) $-isomorphism $%
H_{dR}^{1}(C_{/\mathbb{F}_{q^{2}}})_{-k}\rightarrow \dfrac{V_{k+(q-1)}}{%
D\left( V_{k}\right) }\otimes _{\mathbb{F}}\mathbb{F}_{q^{2}}$ (see proof of
Proposition \ref{long}). Furthermore, the above definition can also be
deduced from Proposition 4.10. in \cite{HaaJan}.

We deduce:

\begin{proposition}
\label{sonno2}For any $q>2$ and any $2\leq k\leq p-1$ we have an isomorphism
of exact sequences of $\mathbb{F}\left[ G\right] $-modules:%
\begin{equation*}
\begin{array}{ccccccccc}
0 & \rightarrow & H^{0}(C_{/\mathbb{F}},\Omega _{C_{/\mathbb{F}}}^{1})_{-k}
& \overset{\iota }{\rightarrow } & H_{dR}^{1}(C_{/\mathbb{F}})_{-k} & 
\overset{\pi }{\rightarrow } & H^{1}(C_{/\mathbb{F}},\mathcal{O}_{C_{/%
\mathbb{F}}})_{-k} & \rightarrow & 0 \\ 
&  & \downarrow \simeq & \circlearrowleft & \downarrow \simeq & 
\circlearrowleft & \downarrow \simeq &  &  \\ 
0 & \rightarrow & e\otimes V_{k-2} & \overset{\overline{\theta }_{q}}{%
\rightarrow } & \dfrac{V_{k+(q-1)}}{D\left( V_{k}\right) } & \overset{\pi
^{\prime \prime }}{\rightarrow } & \limfunc{coker}\overline{\theta }_{q} & 
\rightarrow & 0.%
\end{array}%
\end{equation*}

\noindent Furthermore, if $q=p$ we have $\limfunc{coker}\overline{\theta }%
_{q}\simeq e^{k}\otimes V_{q-1-k}.$
\end{proposition}

\subsection{Integral models}

For a finite field $\mathbb{F}^{\prime }$, let us denote by $W\left( \mathbb{%
F}^{\prime }\right) $ the ring of Witt vectors of $\mathbb{F}^{\prime }$. We
fix a rational prime $l\neq p$, field isomorphisms $\overline{%
\mathbb{Q}
}_{l}\simeq 
\mathbb{C}
$ and $\overline{%
\mathbb{Q}
}_{p}\simeq 
\mathbb{C}
$, and ring embeddings $W\left( \mathbb{F}\right) \hookrightarrow W(\mathbb{F%
}_{q^{2}})\hookrightarrow \overline{%
\mathbb{Q}
}_{p}\simeq 
\mathbb{C}
$. We also denote by $\chi :\mathbb{F}_{q^{2}}^{\times }\rightarrow 
\overline{%
\mathbb{Q}
}_{p}^{\times }$ the Teichm\"{u}ller lifting we considered at the beginning
of section 2, and we let $\Theta \left( \chi ^{k}\right) $ be the cuspidal $%
\overline{%
\mathbb{Q}
}_{p}$-representation of $G=GL_{2}\left( \mathbb{F}\right) $ associated to
the indecomposable character $\chi ^{k}$ ($2\leq k\leq p-1$). We can now go
back to formula \ref{key!}\ and prove the following:

\begin{theorem}
\label{fame2}Let $q>2$ and $2\leq k\leq p-1$ with $k\neq \frac{q+1}{2}$;
there exists an integral model $\widetilde{\Theta }\left( \chi ^{k}\right) $
of the $\overline{%
\mathbb{Q}
}_{p}\left[ G\right] $-module $\Theta \left( \chi ^{k}\right) $ coming from
the $-k$-eigenspace of the first crystalline cohomology group of the
projective curve $C:XY^{q}-X^{q}Y-Z^{q+1}=0$ such that there is an
isomorphism of $\mathbb{F}\left[ G\right] $-modules:%
\begin{equation*}
\dfrac{V_{k+(q-1)}}{D\left( V_{k}\right) }\simeq \overline{\widetilde{\Theta 
}\left( \chi ^{k}\right) }.
\end{equation*}
\end{theorem}

\begin{proof}
As observed in \cite{HaaJan}, \S 2.10, we have for any $g\in U_{2}\left( 
\mathbb{F}_{q^{2}}\right) $:%
\begin{equation*}
\sum\nolimits_{i\geq 0}\left( -1\right) ^{i}\limfunc{tr}\left(
g,H_{cris}^{i}(C_{/\mathbb{F}_{q^{2}}})\otimes _{W(\mathbb{F}_{q^{2}})}%
\overline{%
\mathbb{Q}
}_{p}\right) =\sum\nolimits_{i\geq 0}\left( -1\right) ^{i}\limfunc{tr}\left(
g,H^{i}(C,\overline{%
\mathbb{Q}
}_{l})\right) ,
\end{equation*}

\noindent where $H^{i}\left( \cdot ,\overline{%
\mathbb{Q}
}_{l}\right) $ denotes \'{e}tale cohomology over $\overline{%
\mathbb{Q}
}_{l}$ (for more details, see \cite{DL}). This gives an isomorphisms of $%
U_{2}\left( \mathbb{F}_{q^{2}}\right) $-modules:%
\begin{equation*}
H_{dR}^{1}(C_{/W(\mathbb{F}_{q^{2}})})\otimes _{W(\mathbb{F}_{q^{2}})}%
\mathbb{C}
\simeq H_{cris}^{1}(C_{/\mathbb{F}_{q^{2}}})\otimes _{W(\mathbb{F}_{q^{2}})}%
\mathbb{C}
\simeq H^{1}\left( C,\overline{%
\mathbb{Q}
}_{l}\right) \text{.}
\end{equation*}

\noindent If $1\leq k\leq q$ we deduce $H_{dR}^{1}(C_{/W(\mathbb{F}%
_{q^{2}})})_{-k}\otimes _{W(\mathbb{F}_{q^{2}})}%
\mathbb{C}
\simeq H^{1}\left( C,\overline{%
\mathbb{Q}
}_{l}\right) _{-k}$, hence also an isomorphism of $SL_{2}\left( \mathbb{F}%
\right) $-modules:%
\begin{equation*}
H_{dR}^{1}(C_{/W(\mathbb{F})})_{-k}\otimes _{W(\mathbb{F})}%
\mathbb{C}
\simeq H^{1}\left( C,\overline{%
\mathbb{Q}
}_{l}\right) _{-k}.
\end{equation*}

Since $H^{1}\left( C,\overline{%
\mathbb{Q}
}_{l}\right) _{-k}$ is the subspace of $H^{1}\left( C,\overline{%
\mathbb{Q}
}_{l}\right) $ on which $\mu $ acts via the character $\vartheta _{-k}:\mu
\rightarrow \overline{%
\mathbb{Q}
}_{l}^{\times }:t\mapsto t^{-k}$, by \cite{Lu}, Example 2.20, the $\overline{%
\mathbb{Q}
}_{l}$-representation of $SL_{2}\left( \mathbb{F}\right) $ afforded by $%
H^{1}\left( C,\overline{%
\mathbb{Q}
}_{l}\right) _{-k}$ is of the following type: if $\vartheta _{-k}$ is in
general position (i.e. $\vartheta _{-k}^{2}\neq 1$ or equivalently $k\neq
\left( q+1\right) /2$) $H^{1}\left( C,\overline{%
\mathbb{Q}
}_{l}\right) _{-k}$ is an (irreducible) cuspidal representation of $%
SL_{2}\left( \mathbb{F}\right) $ over $\overline{%
\mathbb{Q}
}_{l}$. If $k=(q+1)/2$, then $H^{1}\left( C,\overline{%
\mathbb{Q}
}_{l}\right) _{-k}=V\oplus V^{\ast }$, with $V$ cuspidal.

From the character theory of $SL_{2}\left( \mathbb{F}\right) $ and $%
GL_{2}\left( \mathbb{F}\right) $ over an algebraically closed field of
characteristic zero (cf. \cite{FH}), we know that, if $\vartheta _{-k}$ is
in general position, there is an indecomposable character $\varsigma :%
\mathbb{F}_{q^{2}}^{\times }\rightarrow 
\mathbb{C}
^{\times }\simeq \overline{%
\mathbb{Q}
}_{l}^{\times }$ for which there is an isomorphism of $SL_{2}\left( \mathbb{F%
}\right) $-modules:%
\begin{equation}
H^{1}\left( C,\overline{%
\mathbb{Q}
}_{l}\right) _{-k}\simeq \limfunc{Res}\nolimits_{SL_{2}\left( \mathbb{F}%
\right) }^{GL_{2}(\mathbb{F})}\left( \Theta \left( \varsigma \right) \right) 
\text{.}  \label{ancora}
\end{equation}

\noindent (Notice that $\varsigma $ is not unique and can be changed by $%
\varsigma ^{-1}$ or any other indecomposable character that equals $%
\varsigma $ on $\mu $). If now we take \textit{any} $p$-adic integral model
of each side of the above isomorphism (e.g. we can take $H_{dR}^{1}(C_{/W(%
\mathbb{F})})_{-k}$ for the \'{e}tale cohomology group) and we reduce $\func{%
mod}p$, we find the $SL_{2}\left( \mathbb{F}\right) $-module isomorphism $%
H_{dR}^{1}(C_{/\mathbb{F}})_{-k}^{ss}\simeq \left( \overline{\widehat{\Theta 
}\left( \varsigma \right) }\right) ^{ss}$ (notice that by Remark \ref%
{reduction} we have a natural identification $H_{dR}^{1}(C_{/\mathbb{F}%
})_{-k}\simeq H_{dR}^{1}\left( C_{/W(\mathbb{F})}\right) _{-k}\otimes _{W}%
\mathbb{F}$). If furthermore $2\leq k\leq p-1$, the first module is
isomorphic to $\left( \frac{V_{k+(q-1)}}{D\left( V_{k}\right) }\right) ^{ss}$
and hence to $\left( \overline{\Theta ^{\prime }\left( \chi ^{k}\right) }%
\right) ^{ss}$ for any choice of integral model $\Theta ^{\prime }\left(
\chi ^{k}\right) $ of $\Theta \left( \chi ^{k}\right) $,\ by formula \ref%
{key!}. Therefore the (reductions of the) Brauer characters of $\overline{%
\widehat{\Theta }\left( \varsigma \right) }$ and $\overline{\Theta ^{\prime
}\left( \chi ^{k}\right) }$ need to coincide.

In the notation of \S 2.2., let $\iota :\mathbb{F}_{q^{2}}\mathbb{%
\hookrightarrow }M_{2}\left( \mathbb{F}\right) $ be given by $c=x+y\sqrt{%
\varepsilon }\mapsto 
\begin{pmatrix}
x & y\varepsilon \\ 
y & x%
\end{pmatrix}%
$, where $x,y\in \mathbb{F}$ and $\varepsilon $ is a generator of $\mathbb{F}%
^{\times }$ (recall $p>2$). If $\iota (c)\in SL_{2}\left( \mathbb{F}\right) $
we have $c\in \mu $; the formulae giving the Brauer characters of the
cuspidal representations of $GL_{2}\left( \mathbb{F}\right) $ imply that, if 
$\varsigma _{|\mu }=\chi _{|_{\mu }}^{h}$ ($0\leq h\leq q$), we have $\chi
(c)^{k}+\chi (c)^{-k}=\chi (c)^{h}+\chi (c)^{-h}$ for any $c\in \mu $, so
that $\left( \chi (c)^{k+h}-1\right) \left( \chi (c)^{k}-\chi (c)^{h}\right)
=0$. We conclude that $k\equiv \pm h(\func{mod}q+1)$ and $\varsigma _{|\mu
}=\chi _{|\mu }^{\pm k}$. We can assume without loss of generality $%
\varsigma =\chi ^{k}$; this implies that the $SL_{2}\left( \mathbb{F}\right) 
$-action on $H^{1}\left( C,\overline{%
\mathbb{Q}
}_{l}\right) _{-k}$ extends to a $GL_{2}(\mathbb{F})$-action giving an
isomorphism $H^{1}\left( C,\overline{%
\mathbb{Q}
}_{l}\right) _{-k}\simeq \Theta \left( \chi ^{k}\right) $.

If $\widetilde{\Theta }\left( \chi ^{k}\right) $ is the $W\left( \mathbb{F}%
\right) $-model of $\Theta \left( \chi ^{k}\right) $ corresponding to $%
H_{dR}^{1}(C_{/W(\mathbb{F)}})_{-k}$ in the above isomorphism, we have $%
\overline{\widetilde{\Theta }\left( \chi ^{k}\right) }\simeq H_{dR}^{1}(C_{/%
\mathbb{F}})_{-k}$.
\end{proof}

\section{Modular forms}

In this section we assume that $q=p$ is a prime number larger than $3$ and,
for any non-negative integer $k$, we denote by $V_{k}$ the $\overline{%
\mathbb{F}}_{p}[GL_{2}(\mathbb{F}_{p})]$-module $\limfunc{Sym}\nolimits^{k}%
\overline{\mathbb{F}}_{p}^{2}$.

Let us fix an integer $N\geq 5$ not divisible by the prime $p$; let $%
\mathcal{E}\rightarrow X_{1}\left( N\right) $ be the universal generalized
elliptic curve over the modular curve $X_{1}\left( N\right) $ (all the
schemes here are over $\overline{\mathbb{F}}_{p}$) and let $\omega
_{X_{1}(N)}=0^{\ast }\Omega _{\mathcal{E}/X_{1}\left( N\right) }^{1}$ be the
coherent sheaf obtained by pulling-back via the zero section $0:X_{1}\left(
N\right) \rightarrow \mathcal{E}$ the sheaf of K\"{a}hler differentials of $%
\mathcal{E}$ over $X_{1}\left( N\right) $. For any integer $k\geq 0$ we let $%
M_{k}\left( \Gamma ,\overline{\mathbb{F}}_{p}\right) :=H^{0}(X_{1}\left(
N\right) ,\omega _{X_{1}(N)}^{\otimes k})$ be the space of $\func{mod}p$
modular forms of weight $k$ and level $\Gamma :=\Gamma _{1}(N)$.

In view of the Eichler-Shimura isomorphism and of \cite{A-S}, Proposition
2.5., the study of Hecke eigensystems of $\func{mod}p$ modular forms of
weight $k\geq 2$ and level $N$ leads naturally to the study of the
eigenvalues of the Hecke algebra $\mathcal{H}_{N}$ acting on the cohomology
group $H^{1}(\Gamma ,V_{k-2})$, where $\Gamma $ acts on $V_{k-2}$ via its
reduction $\func{mod}p$, and the action of $\mathcal{H}_{N}$ is defined as
in \cite{A-S} (here $\mathcal{H}_{N}$ is generated over $\overline{\mathbb{F}%
}_{p}$ by the Hecke operators $T_{r}$ where $p\nmid r$).

\bigskip

As for the group $GL_{2}\left( \mathbb{F}_{p}\right) $ we have two operators
raising the levels of the symmetric power representations $V_{k}$'s by $p+1$
and $p-1$ respectively, the theory of modular forms $\func{mod}p$ offers two
maps that increase weights by $p+1$ and $p-1$: the theta operator $\Theta $
and the Hasse invariant $A$ respectively (cf. \cite{Se3} for their
definitions). One might expect a connection between the representation
theoretical side discussed in the previous sections and the modular forms
side.

Indeed, in \cite{A-S}, A. Ash and G. Stevens identified a group-theoretical
analogue of the $\Theta $-operator (\cite{A-S}, Theorem 3.4.) in the map
induced in cohomology by the operator $\theta _{p}$, that is by the
multiplication by the Dickson polynomial $X^{p}Y-XY^{p}$ :\ for any $k\geq 2$
the map $\theta _{p}$ induces an Hecke-equivariant map:%
\begin{equation*}
\theta _{p,\ast }:H^{1}(\Gamma ,V_{k-2})\rightarrow H^{1}(\Gamma ,V_{k+p-1})
\end{equation*}

\noindent that corresponds, on the space of modular forms, to raising the
weight by $p+1$.

\bigskip

\subsection{Group cohomology and the Hasse invariant}

In \cite{K-E}, Edixhoven and Khare construct a cohomological analogue of the
Hasse invariant: by studying the degeneracy map $H^{1}\left( \Gamma
,V_{0}\right) ^{2}\rightarrow H^{1}\left( \Gamma \cap \Gamma _{0}\left(
p\right) ,V_{p-1}\right) $, they determine a monic Hecke-equivariant
homomorphism:%
\begin{equation*}
\alpha :H^{1}(\Gamma ,V_{0})\hookrightarrow H^{1}(\Gamma ,V_{p-1}).
\end{equation*}%
Notice that $\alpha $ is not defined in \cite{K-E} as coming from a morphism
on coefficients.

The $D$-map allows us to extend the existence of an injection as the map $%
\alpha $ above:

\begin{proposition}
\label{ora!}Let $\mathfrak{M}$ be a non-Eisenstein maximal ideal of $%
\mathcal{H}_{N}$:

\begin{enumerate}
\item if $k\geq 0$ and $H^{1}(\Gamma ,V_{k})_{\mathfrak{M}}\neq 0$, then
also $H^{1}(\Gamma ,V_{k+(p-1)})_{\mathfrak{M}}\neq 0$;

\item if $0\leq k\leq p-1$, there is an embedding of Hecke modules:%
\begin{equation*}
H^{1}(\Gamma ,V_{k})_{\mathfrak{M}}\hookrightarrow H^{1}(\Gamma
,V_{k+(p-1)})_{\mathfrak{M}}
\end{equation*}

\noindent \noindent that is induced by $D$ if $0<k\leq p-1$, and is the
above map $\alpha $ for $k=0$.
\end{enumerate}
\end{proposition}

\begin{proof}
If $k\geq 0$ and $k\neq 0(\func{mod}p+1)$, Proposition \ref{negativity}
applied with $q=p$ (and a suitable re-indexing) says that $V_{k+(p-1)}-V_{k}$
is positive in $K_{0}\left( GL_{2}(\mathbb{F}_{p})\right) $, giving the
first assertion for $k\neq 0(\func{mod}p+1)$.

If $1\leq k\leq p-1$ we have the exact sequence of $GL_{2}(\mathbb{F}_{p})$%
-modules:%
\begin{equation*}
0\rightarrow V_{k}\overset{D}{\rightarrow }V_{k+(p-1)}\rightarrow \limfunc{%
coker}D\rightarrow 0\text{.}
\end{equation*}

\noindent By passing to the long exact sequence in cohomology and localizing
with respect to the non-Eisenstein maximal ideal $\mathfrak{M}$ we get the
second statement for $1\leq k\leq p-1$ (cf. \cite{K}). If $k=0$, the
existence of a monic map $\alpha :H^{1}(\Gamma ,\mathbb{F}_{p})_{\mathfrak{M}%
}\hookrightarrow H^{1}(\Gamma ,V_{p-1})_{\mathfrak{M}}$ is the cited above
result of Edixhoven and Khare (\cite{K-E}).

The existence of $\alpha $ also implies the first statement for $k\equiv 0(%
\func{mod}p+1)$: if $k=s(p+1)$ for some $s\geq 0$, formula \ref{questa!!!}
gives the following identity in $K_{0}\left( GL_{2}(\mathbb{F}_{p})\right) $:%
\begin{equation}
V_{s\left( p+1\right) +p-1}=e^{s}\cdot V_{p-1}+(V_{s(p+1)}-e^{s}\cdot V_{0}).
\label{vaii}
\end{equation}

Notice that $V_{s(p+1)}-e^{s}\cdot V_{0}>0$ because of the existence of the
monic map $\theta _{p}:e\otimes V_{0}\hookrightarrow V_{p+1}$. If $%
H^{1}(\Gamma ,V_{s(p+1)})_{\mathfrak{M}}\neq 0$ then $H^{1}(\Gamma
,e^{s}\otimes V_{0})_{\mathfrak{M}}\neq 0$ or $H^{1}(\Gamma
,V_{s(p+1)}/e^{s}\otimes V_{0})_{\mathfrak{M}}\neq 0$; in the first case, by
applying $\alpha $ to the twisted module $H^{1}(\Gamma ,e^{s}\otimes V_{0})_{%
\mathfrak{M}}$, we deduce $H^{1}(\Gamma ,e^{s}\otimes V_{p-1})_{\mathfrak{M}%
}\neq 0$ and hence, by \ref{vaii} above, $H^{1}(\Gamma ,V_{s(p+1)+p-1})_{%
\mathfrak{M}}\neq 0$. If it is $H^{1}(\Gamma ,V_{s(p+1)}/e^{s}\otimes
V_{0})_{\mathfrak{M}}\neq 0$, we conclude again $H^{1}(\Gamma
,V_{s(p+1)+p-1})_{\mathfrak{M}}\neq 0$.
\end{proof}

\begin{remark}
The above result cannot be deduced only by the existence of the map $D$ and
by Proposition \ref{negativity}, since if $k=0$ then $V_{p-1}-V_{0}$ is not
positive in $K_{0}\left( GL_{2}(\mathbb{F}_{p})\right) $.
\end{remark}

\bigskip

\subsection{\textbf{Construction of the maps }$\protect\alpha _{k}$\textbf{\
for }$k>0$}

The definition of $\alpha $ is given in \cite{K-E} for $k=0$; however
starting from this case it is not hard to generalize the construction to any
integer $k\geq 0$. Recall that we fixed a prime $p>3$ and an integer $N\geq
5 $ prime to $p$.

Let $\Gamma _{0}:=\Gamma \cap \Gamma _{0}(p)$, $g:=%
\begin{pmatrix}
p & 0 \\ 
0 & 1%
\end{pmatrix}%
$ and $g^{\ast }=%
\begin{pmatrix}
1 & 0 \\ 
0 & p%
\end{pmatrix}%
$; denote by $\beta ^{\prime }:H^{1}(\Gamma ,V_{k})\rightarrow H^{1}(\Gamma
_{0},V_{k})$ the restriction map twisted by $g$, i.e. the homomorphism
induced by assigning to the 1-cocycle $\zeta \in Z^{1}(\Gamma ,V_{k})$ the
cocycle: 
\begin{equation*}
x\mapsto g^{\ast }\cdot \zeta (gxg^{-1}),\text{ \ \ \ }x\in \Gamma _{0}.
\end{equation*}%
(For $k=0$ the action of $g^{\ast }$ is trivial, and $\beta ^{\prime }$ is
the twisted restriction considered in \cite{K-E}; if $k>0$, $g^{\ast }$ acts
as multiplication by the elementary matrix $E_{11}\in M_{2}(\mathbb{F}_{p})$%
). Notice that by Shapiro's lemma, we have a natural isomorphism of Hecke
modules $\beta ^{\prime \prime }:H^{1}(\Gamma _{0},V_{k})\rightarrow
H^{1}(\Gamma ,\limfunc{Ind}_{\Gamma _{0}}^{\Gamma }V_{k})$.

Let us now fix a left transversal $\left\{ x_{1},...,x_{p+1}\right\} $ of $%
\Gamma _{0}$ inside $\Gamma $; e.g. we can take:%
\begin{equation}
x_{i}:=%
\begin{pmatrix}
1 & 0 \\ 
N(i-1) & 1%
\end{pmatrix}%
\text{ (}1\leq i\leq p\text{), \ }x_{p+1}=%
\begin{pmatrix}
Bp & -1 \\ 
AN & 1%
\end{pmatrix}%
\text{, }  \label{transversal}
\end{equation}

\noindent where $AN+Bp=1$ and $A,B\in 
\mathbb{Z}
$.\ Denote by $\mathbb{F}_{p}\left[ \Gamma /\Gamma _{0}\right] $ the $%
\mathbb{F}_{p}$-vector space of functions $f:\Gamma /\Gamma _{0}\rightarrow 
\mathbb{F}_{p}$ endowed with the left $\Gamma $-action defined by $\left(
x\cdot f\right) \left( x_{i}\Gamma _{0}\right) =f\left( x^{-1}x_{i}\Gamma
_{0}\right) $, for every $x\in \Gamma ,f\in \mathbb{F}_{p}\left[ \Gamma
/\Gamma _{0}\right] $ and $1\leq i\leq p+1$. There is a $\Gamma $%
-isomorphism $b:\limfunc{Ind}_{\Gamma _{0}}^{\Gamma }V_{k}\rightarrow 
\mathbb{F}_{p}\left[ \Gamma /\Gamma _{0}\right] \otimes _{\mathbb{F}%
_{p}}V_{k}$ (where the codomain is endowed with the diagonal $\Gamma $%
-action) given by sending an element $f:\Gamma /\Gamma _{0}\rightarrow V_{k}$
of $\limfunc{Ind}_{\Gamma _{0}}^{\Gamma }V_{k}$ to:%
\begin{equation*}
\sum\nolimits_{i=1}^{p+1}x_{i}\Gamma _{0}\otimes x_{i}f\left( x_{i}\Gamma
_{0}\right) .
\end{equation*}

Fixing the choice of transversal given by \ref{transversal},\ we can
identify $\mathbb{F}_{p}\left[ \Gamma /\Gamma _{0}\right] $ and $\mathbb{F}%
_{p}\left[ \mathbb{P}^{1}\left( \mathbb{F}_{p}\right) \right] $ through the
bijection $\Gamma /\Gamma _{0}\rightarrow \mathbb{P}^{1}\left( \mathbb{F}%
_{p}\right) $ defined by $x_{i}\mapsto (1:N(i-1))$ for $1\leq i\leq p$, and $%
x_{p+1}\longmapsto (0:1)$. Using Lemma \ref{questooo}, we find: 
\begin{equation*}
\mathbb{F}_{p}\left[ \Gamma /\Gamma _{0}\right] \otimes _{\mathbb{F}%
_{p}}V_{k}\simeq \mathbb{F}_{p}\left[ \mathbb{P}^{1}\left( \mathbb{F}%
_{p}\right) \right] \otimes _{\mathbb{F}_{p}}V_{k}\simeq V_{k}\oplus
(V_{p-1}\otimes _{\mathbb{F}_{p}}V_{k}).
\end{equation*}

\noindent The composition of $\beta ^{\prime \prime }\circ \beta ^{\prime }$
with the map induced in cohomology by $b$, and the above isomorphisms yield
an Hecke equivariant map: 
\begin{equation*}
\beta :H^{1}(\Gamma ,V_{k})\rightarrow H^{1}(\Gamma ,V_{k})\oplus
H^{1}(\Gamma ,V_{p-1}\otimes _{\mathbb{F}_{p}}V_{k}).
\end{equation*}

\noindent By composing $\beta $ with the projection onto the second factor
and then with the map induced in cohomology by polynomial multiplication $%
V_{p-1}\otimes _{\mathbb{F}_{p}}V_{k}\rightarrow V_{k+(p-1)}$, we finally
obtain:

\begin{proposition}
For every integer $k\geq 0$ there is an Hecke-equivariant homomorphism:%
\begin{equation*}
\alpha _{k}:H^{1}(\Gamma ,V_{k})\rightarrow H^{1}(\Gamma ,V_{k+(p-1)})
\end{equation*}

\noindent such that $\alpha _{0}=\alpha $ is the map defined in \cite{K-E}.
\end{proposition}

Some questions arise naturally. For example, we know that $\alpha _{0}$ is
monic, and one might expect that also $\alpha _{k}$ is monic in the range $%
1\leq k\leq p-1$, i.e. when $V_{k}\neq \mathbb{F}_{p}$ is irreducible. In
this range, by Proposition \ref{ora!}, we also have - for any non-Eisenstein
maximal ideal $\mathfrak{M}$ of $\mathcal{H}_{N}$ - a monic Hecke-map $%
H^{1}(\Gamma ,V_{k})_{\mathfrak{M}}\hookrightarrow H^{1}(\Gamma
,V_{k+(p-1)})_{\mathfrak{M}}$ induced by $D$, that one would like to compare
with $\alpha _{k}$.

To this purpose, it could be useful to see the maps $\alpha _{k}$'s as
homomorphisms between groups of modular symbols, using the Hecke-equivariant
isomorphism%
\begin{equation*}
\limfunc{Hom}\nolimits_{\Gamma }(\mathcal{D}_{0},V_{k})\overset{\simeq }{%
\rightarrow }H_{c}^{1}(\Gamma ,V_{k}),
\end{equation*}

\noindent where $\mathcal{D}_{0}$ is the group of divisors of degree zero of 
$\mathbb{P}^{1}\left( 
\mathbb{Q}
\right) $, and $H_{c}^{1}(\Gamma ,\cdot )$ denotes compactly supported
cohomology (cf. \cite{A-S}, \S 4).

\end{document}